\documentclass[preprint,12pt]{elsarticle}

\usepackage{amssymb}
\usepackage{amsmath}
\usepackage{amsthm}
\usepackage{mathtools}
\usepackage{dsfont}
\usepackage[capitalise]{cleveref}
\usepackage[printonlyused,withpage]{acronym}
\usepackage{subcaption,graphicx}
\usepackage{epstopdf}
\usepackage{algorithm}
\usepackage{algorithmicx}
\usepackage{algpseudocode}
\usepackage{multirow}

\usepackage{nomencl}
\makenomenclature

\newcommand{\up}[1]{$^{\mathrm{#1}}$}
\newcommand{\tol}{\mathsf{tol}}
\newcommand{\y}{\mathrm{y}}
\newcommand{\hy}{\mathrm{y}^*}
\newcommand{\e}{{e}}
\newcommand{\h}{\tau}
\renewcommand{\t}{\mathrm{t}}
\newcommand{\I}{\mathrm{I}}
\newcommand{\U}{\mathcal{U}}
\newcommand{\V}{\mathcal{V}}
\newcommand{\J}{\mathcal{J}}
\newcommand{\f}{{f}}
\newcommand{\C}[1]{\mathrm{C}_{#1}}
\newcommand{\T}[1]{\mathbf{T}_{#1}}
\newcommand{\Cc}[1]{\mathcal{C}_{#1}}
\newcommand{\F}[2]{\mathrm{F}_{#1}^{#2}}
\newcommand{\Mp}{\boldsymbol{\mathrm{mp}}}
\newcommand{\FE}{\boldsymbol{\mathrm{RK1}}}
\newcommand{\RK}[1]{\boldsymbol{\mathrm{RK}#1}}
\newcommand{\CRK}[1]{\boldsymbol{\mathrm{CRK}#1}}

\newcommand{\CGRK}[1]{\boldsymbol{\mathrm{CGRK}#1}}
\newcommand{\DOP}[1]{\boldsymbol{\mathrm{DOPR}#1}}
\newcommand{\BS}[1]{\boldsymbol{\mathrm{BS}#1}}
\newcommand{\BPL}{\boldsymbol{\mathrm{BPL}}}
\newcommand{\pert}{\varepsilon}
\newcommand{\Op}{\mathcal{O}}
\renewcommand{\Re}{\mathrm{Re}}
\renewcommand{\Im}{\mathrm{Im}}
\renewcommand{\u}{\mathrm{u}}
\renewcommand{\i}{\mathrm{i}}
\renewcommand{\v}{\mathrm{v}}
\renewcommand{\d}{\mathrm{d}}
\newcommand{\op}{\mathrm{p}}
\renewcommand{\S}[1]{\mathrm{S}_{#1}}
\renewcommand{\r}{\mathrm{r}}
\renewcommand{\a}{\mathrm{a}}
\renewcommand{\b}{\mathrm{b}}
\renewcommand{\c}{\mathrm{c}}
\newcommand{\g}{\mathrm{g}}
\newcommand{\w}{\mathrm{w}}
\renewcommand{\H}{\mathrm{H}}
\newcommand{\D}{\mathcal{D}}
\newcommand{\Dc}{\mathcal{D}_c}
\renewcommand{\P}{\mathrm{P}}
\newcommand{\Pc}{\mathrm{P}_c}

\newcommand\dg\mathfrak
\newcommand{\B}{\mathcal{B}}

\newcommand{\Res}{\mathcal{R}es}

\newtheorem{theorem}{Theorem}
\newtheorem{lemma}[theorem]{Lemma}
\newdefinition{example}{Example}

\begin{document}

\begin{frontmatter}



\title{Error estimation for numerical approximations of ODEs via composition techniques. Part I: One-step methods}

\author[label1]{Ahmad Deeb}
\ead{ahmad.deeb@ku.ac.ae}
\author[label1,label2]{Denys Dutykh}
\ead{denys.dutykh@ku.ac.ae}
\affiliation[label1]{organization={Khalifa University of Science and Technology},
            addressline={PO Box 127788},
            city={Abu Dhabi},
            country={UAE}}

\affiliation[label2]{organization={Causal Dynamics},
            addressline={Pty LTD},
            city={Perth},
            country={Australia}}

\begin{abstract}
In this study, we introduce a refined method for ascertaining error estimations in numerical simulations of dynamical systems via an innovative application of composition techniques. Our approach involves a dual application of a basic one-step numerical method of order $\op$ in this part, and for the class of Backward Difference Formulas schemes in the second part [Deeb A., Dutykh D. and AL Zohbi M. \textit{Error estimation for numerical approximations of ODEs via composition techniques. Part II: BDF methods}, Submitted, 2024]. This dual application uses complex coefficients, resulting outputs in the complex plane. The methods innovation lies in the demonstration that the real parts of these outputs correspond to approximations of the solutions with an enhanced order of $\op+1$, while the imaginary parts serve as error estimations of the same order, a novel proof presented herein using Taylor expansion and perturbation technique. The linear stability of the resulted scheme is enhanced compared to the basic one. The performance of the composition in computing the approximation is also compared. Results show that the proposed technique provide higher accuracy with less computational time. This dual composition technique has been rigorously applied to a variety of dynamical problems, showcasing its efficacy in adapting the time step, particularly in situations where numerical schemes do not have theoretical error estimation. Consequently, the technique holds potential for advancing adaptive time-stepping strategies in numerical simulations, an area where accurate local error estimation is crucial yet often challenging to obtain.\end{abstract}




\end{frontmatter}


\section{Introduction}
\label{sec1}

Solving \ac{ODE}s and \ac{PDE}s by numerical methods is a powerful tool when there is no possibility to have the exact solution \cite{book:butcher,book:tomas}. Numerical schemes were designed for the integration of both stiff and non-stiff problems \cite{book:hairer}. To achieve high accuracy and improve stability of numerical solutions \cite{wang-2003}, it is preferred to use ''\textit{variable}'' or ''\textit{adaptive}'' time stepping techniques to follow the dynamics of solutions \cite{book:hairer,book:iserles,DEEB_2022_bpl}. Thus, the time step should be decreased when the system presents high stiffness and increased when dynamics is locally simple.
Stiff problems could be defined when explicit methods fail for a high stiffness ratio while implicit methods succeed \cite{curtiss_1952}.
Different methods exist for integrating stiff problems such as the \ac{GRK} methods or the \ac{BDF} \cite{book:hairer2}. \ac{ATS} could be also achieved using, if possible, the numerical energy stability \cite{ZHAO2020124901,HUANG2023127622}, the dynamics of the residual error or the rate of change of the solution \cite{RC-CN}.

Another way to adapt the time step is by having a local error estimate of the approximation \cite{zennaro-85,wang-2003} produced by a given numerical scheme. The time step is chosen to meet a user defined tolerance $\tol$.
The first who wrote about changing the time step size during simulation to accelerate computation was Runge \cite{Runge_1895}. \ac{ERK} methods are nowadays the most useful methods for \ac{ATS}. These methods consist of adding coefficients $\{b_i^*\}$ in their Butcher tables for predicting additional lower order approximation $\hy_n$ without requiring additional function evaluations.
The difference $\e^*_n = \|\y_n-\hy_n\|$, where $\y_n$ is the first approximation using coefficients $b_i$, is considered as an error estimate of $\e(\t_n) =
\y(t_n)-\hy_n$. Having already computed $K_i$ to evaluate the first approximation $\y_n$, $\hy_n$  is obtained via arithmetic operations involving the terms $K_i$,  $b_i^*$ and the time step $\h$. We list here, for instance, the Bogacki-Shampine method \cite{bogacki_shampine} that computes the second and third-order approximations using four stages. The Dormand-Prince \cite{dromand_prince} method computes the fourth and the fifth-order accurate solutions to produce a fourth order error estimate. Alternative methods for error estimation were developed in the context of parabolic problems \cite{johnson-90}, linear hyperbolic equations \cite{xu-20}, linear Schrödinger equation \cite{Liao-10}, Maxwell's equations \cite{ZHOU201832} and Maxwell-Schrodinger equations \cite{ma-18}. The error estimates are used in the \ac{ATS} technique for solving problems in fluid mechanics \cite{ranocha-22}, in thermal mechanics \cite{Monge-2020}, wave propagation \cite{Dutykh-2015}, and electromagnetic \cite{Dinavahi-2021,Dinavahi-17}, where the new time step is generally evaluated as follows:
\begin{equation}
 \label{adapting-step}
 \h_{n+1} = C \times \h_n\left(\frac{\tol}{\e_n}\right)^{\frac{1}{\op+1}}.
\end{equation}
Here $C<1$ is a safety factor and $\op$ is a parameter that is, by default, the order of the numerical scheme. \ac{ERK} methods are not the only ones that provide error estimates for implementing the \ac{ATS}.
It can be also implemented using a priori error estimates of the local error by $\e_n \simeq \C{}\times  \h_n^{\op+1}$, where $\C{}$ is a positive error constant to be determined for every one-step method. \ac{LMS} methods were also studied to provide  a priori error estimations of their approximations.
We refer the reader to \cite{book:hairer, book:hairer2, hairer2002geometric} for additional details.

Back to formula \eqref{adapting-step}, one may ask the question: {\bf could we provide an error estimate for every time marching numerical scheme we use}? What if we conceive a new one-step method and we do not yet have theoretical estimations for its numerical approximations, {\bf is there any tool to predict its error}? The answers to both questions are affirmatives using the difference between two approximations obtained by a higher-order numerical scheme and the original one. Nevertheless, this will be to the detriment of additional computations cost while computing the higher-order approximations.

In this work, we are interested in providing a numerical solution that possesses an error estimate to be used later for the time stepping technique. This is done by using a double jump composition. The composition technique is commonly used to increase the order of approximation of basic one-step methods. Here, we will use it with complex coefficients, \emph{i.e.} the distance between two consecutive iterations is a time instant lying in the complex plane. At the end of the double composing the one-step method, the process produces outputs with values in the complex plane. The real parts are considered to be the new approximations of order $\op+1$ \cite{casas_2021_complex}, while the imaginary parts are proven here to be error estimates, of the real ones, having the accuracy of order $\op+1$ too. In the second part \cite{deeb:part2} of this work, we will provide error estimates by extending the composition technique to a class of \ac{LMS}, namely the \ac{BDF} schemes.

For any one-step method, we associate a function $\Phi_{\h}$ to be its numerical flow. To integrate the solution from $\t$ to $\t+\h$, we apply  $\Phi_{\h}$ on such an approximation at instant $\t$ with a time step $\h$. The composition technique consists of applying $\Phi_{\gamma_i\h}$ successively $s$ times $\{\gamma_i \in \mathds{C}\,\vert\, i\in \S{1}^{s}\}$, where the time step at each sub-layer is equal to $\gamma_i \h$. The objective of this composition was first proposed in \cite{Mc-1995,yoshida-1990,suzuki-1990} to design a new numerical flow having properties of symmetry and symplecticity. For example, the St{\"o}rmer-Verlet method \cite{hairer_lubich_wanner_2003}, which is a symmetric and symplectic scheme used in simulating Hamiltonian systems \cite{leimkuhler_reich_2005}, is the result of composing two schemes with a half time step: the symplectic Euler scheme \cite{cromer-81}, and its adjoint defined by its inverse with a negative time step. The outcome numerical flow was able to produce an approximation with a greater accuracy order than the initial one.
Another composition was done by Butcher \cite{butcher-1969}. He proposed a fifth-order scheme with five stages, by composing two \ac{RK} methods of order three. A generalization was extended to $\mathbf{B}$-series (see \cite[section II.2]{book:hairer} and \cite{iserles-1984}) to overcome stability issues when low-order numerical schemes were used.
This technique was applied to solve numerically dynamical systems as in electromagnetic \cite{channell-14}, in quantum mechanics for Klein-Gordon lattices \cite{sokos-18,bader-19}, in astronomy \cite{blanes-13}, in electrical and chaotic systems \cite{butusov-16}, \emph{etc}.

The works of Blanes \emph{et al.} \cite{casas-2006,casas-2008} developed for a family of coefficients satisfying algebraic conditions the framework of constructing families of higher-order numerical integrators by composing basic ones. In a recent paper, Casas \emph{et al.} \cite{casas_2021_complex} constructed and analyzed a new class of numerical integrators by composing twice a basic one-step integrator with complex coefficients. They showed that this composition produces a new numerical integrator with an additional accuracy order and that the symmetric and symplectic properties are preserved up to higher-orders. The use of complex coefficients are extended to develop symmetric conjugate splitting methods \cite{blanes-24-1} for solving linear parabolic evolution problems, as illustrated by the linear Schr\"odinger equations, where a splitting method was designed when the linear equation is split with two linear operators: $A$ and $B$. It was shown on one hand that the use of complex coefficients and their conjugates in designing symmetric-conjugate splitting method maintains bounded errors over time. On the other hand, they showed that the imaginary part of the approximation presents an error estimation of the solution. However, their illustration is limited to linear problems. In addition, their analysis showing an error estimation in the imaginary part works only when the commutator of the two linear operators, $A$ and $B$, is not null. If the latter is not valid, the imaginary part in the produced approximation is null, thus the approximation of the error could not be retrieved. In the case of non-linear problems, a modified splitting technique was applied in the case of Gross–Pitaevskii systems \cite{blanes-24-2}, with real and complex coefficients. Nonetheless, nothing was mentioned about error estimation using the imaginary part, thus the urgent need in providing theoretical proofs of providing error estimator in the imaginary part for non-linear problems.

In this manuscript (first part of this work), we provide an error estimation of a new numerical approximation obtained by composing twice any one-step method of order $\op$. The use of complex coefficients, verifying algebraic equations, will provide outputs in the complex plane. The real part of the output is considered to be the new approximation of order $\op+1$ that was already shown in \cite{casas_2021_complex}. The novelty in this paper is the proof, by techniques of perturbations, that the imaginary part presents an error estimate of the constructed approximation by any one step method applied to non-linear differential equations. An error constant is provided too. The linear stability of the resulted scheme by the real part of the composition is studied for different types of schemes. The presented examples show enhancing in linear stability exhibited in the increasing of the stability domain in the left part of the complex plane. The performance of the composed scheme, for several examples, is exhibited too and compared to the basic integrator. Despite the use of arithmetic complex by the composition, composed schemes outperform the basic integrators in providing approximations with higher accuracy and lower computational time.

The outline of this paper will be as follows. \cref{sec2} presents the mathematical framework and results proving that the imaginary part is an error estimate of order $\op+1$.
\cref{sec3} will present some of familiar numerical integrators, where their linear stability are compared with the scheme resulted by the real part of their composition. The rate of convergence for both, before and after the composition, are shown. The performance of basic integrators and their composition is presented in this section too.
\cref{sec4} exhibits employing the composition of the above schemes in solving some \ac{ODE}s, showing the efficiency of the proposed technique in CPU and error estimate. We end with main conclusions and some perspectives of this work highlighted in \cref{sec5}.

\section{Mathematical framework and results}
\label{sec2}

Consider the following \ac{CP} to solve:
\begin{equation}
 \label{diff-system}
 \frac{\d\y}{\d\t} = \f(\t,\y),\qquad \y(\t_0) = \y_0,
\end{equation}
with
\begin{equation}
\y: \left\lvert
 \begin{array}{ccc}
  \I\subset \mathds{R} &\longrightarrow &\U \subseteq \mathds{R}^d,\\
  \t &\longmapsto &\y(\t),
 \end{array}
 \right.
\qquad
\f: \left\lvert
\begin{array}{ccc}
  \I \times \U&\longrightarrow &\mathds{R}^d,\\
  (\t,\y) &\longmapsto &\f(\t,\y(t)).
 \end{array}
 \right.
 \end{equation}
We consider the cases where $\f$ could be naturally extended to $\J\times \V$, where $\I\subset \J \subset \mathds{C}$ and $\U\subset \V\subset \mathds{C}^d $.
Exact solutions are in most cases impossible to find, thus numerical approximations are sought on a discrete, ordered set of points $\{\t_0,\t_1,\ldots,\t_n,\ldots\} \subset \I$ with $\t_i<\t_{i+1}, \, \forall i=0,1,\ldots$. We denote by $\y_n$ the approximation of $\y(\t_n)$ and by $\h_n = \t_{n+1}-\t_n$ the $n$\up{th} step size. For every $\t \in \I$, we denote by $\varphi_{\t}$ the following flow map:
\begin{equation}
 \varphi_{\t}:
 \left\lvert
 \begin{array}{ccc}
  \U\subset\mathds{R}^d &\longrightarrow &\U\subset\mathds{R}^d,\\
  \y_0 & \longmapsto & \varphi_{\t}(\y_0) = \y(\t),
 \end{array}
 \right.
\end{equation}
having the solution of the above \ac{CP} with a given initial condition $\y_0 \in \mathds{R}^d$ as the image. This map is called the exact flow of the \ac{IVP}. 

Many algorithms and numerical integration schemes for system \eqref{diff-system} were proposed to approximate their solutions \cite{book:hairer, book:hairer2}. Generally, there are two classes of schemes: one-step and \ac{LMS} methods which are not considered in this part. For one-step numerical integration, it predicts the solution at the following instant $\t_{n+1}$ by using only the last known approximation, as for instance the famous explicit fourth-order \ac{RK} scheme \cite{Runge_1895,butcher_1964}. Having a one-step numerical method and a set of equidistant points $\t_n$ ($\h_n \equiv \h, \forall n \in \mathds{N}$), we can associate a numerical flow, denoted by $\Phi_\h$, such that:
\begin{equation}
 \Phi_{\h}:\left\lvert
 \begin{array}{ccc}
  \mathds{R}^d &\longrightarrow &\mathds{R}^d\\
  \y_n & \longmapsto & \Phi_\h(\y_n) = \y_{n+1}
 \end{array}\right.,
\end{equation}
where for every $\t_{n+1}\ \coloneqq\ \t_0+(n+1)\h$, the solution $\y(\t_{n+1})$ is approximated by the image of $\Phi_{\h}\left(\y_{n}\right)$ of the associated numerical flow. We say that the numerical flow $\Phi_\h$ is of order $\op$ if the local error follows the asymptotic relation:
\begin{equation}
 \e(\t_n) \coloneqq \y(\t_n) - \y_n = O(\h^{\op+1}).
\end{equation}
To be more precise, one can write the following relation for every one-step method of order $\op$, by assuming that we have $\y_n \equiv \y(t_n)$:
\begin{equation}
\label{num-flow}
\begin{array}{rcl}
 \Phi_{\h}(\y_n) = \varphi_\h(\y_n) &+& \C{\op+1} \h^{\op+1} \F{\op+1}{}(\t_n,\y_n)\\
 &+& \C{\op+2} \h^{\op+2} \F{\op+2}{}(\t_n,\y_n) + O(\h^{\op+3}),
\end{array}
 \end{equation}
where $\F{\op+1}{}(\t_n)$ is a function of $\y_n$ that it is expressed through $\f(\cdot,\cdot)$ and its derivatives.

\begin{example}
The midpoint scheme has the associated numerical flow $\Phi^{\Mp}_{\h}$:
\begin{equation}
 \Phi^{\Mp}_{\h}(\y_n) = \y_n + \h \f\left( \t_n+\frac{\h}{2}, \y_n +\frac{\h}{2}\f(\t_n , \y_n) \right).
\end{equation}
After writing the Taylor expansion of  the function $\f(\cdot,\cdot)$ in the neighborhood of $\t_n$ and $\y_n$, assembling terms of $\h$ and its powers, we can write the following equality for the midpoint rule:
\begin{eqnarray*}
 \Phi^{\Mp}_{\h}(\y_n) &=& \varphi_{\h} (\y_n) -\frac{1}{24}\h^3 \F{3}{\Mp}(\t_n,\y_n)  \\
 &&-\frac{1}{48}\h^4 \F{4}{\Mp}(\t_n,\y_n) +  O(\h^{5}),\\
  \F{3}{\Mp}(\t,\y) &=& \big(\f^{(2)} + 3\f_{y}\f^{\prime}\big)(\t,\y)\\ 
  \F{4}{\Mp}(\t,\y) &=&\Big(\f^{(3)} +\f_\y\f^{(2)} + 3(\f_{ty} + \f_{yy}\f)\f^{\prime}\Big)(\t,\y).
 \end{eqnarray*}
 \end{example}
We state the following lemma to prepare the proof of our main theorem:
 \begin{lemma}
 \label{lemma1}
  Consider any numerical flow associated with a one-step method of order $\op$ where  \cref{num-flow} follows for any state $\y_n$. If we perturb $\y_n$ by $\pert \sim\Op(\h^{\op+1})$, then, the numerical flow applied to the perturbed element satisfies the relation below:
  \begin{equation}
   \label{pert-num-flow}
   \Phi_{\h}(\y_n + \pert) = \Phi_{\h}(\y_n) + \pert + \pert\, \h \,\f_\y(\t_n,\y_n) + \Op(\h^{\op+3}).
  \end{equation}
 \end{lemma}
 \begin{proof}
  To establish the proof, we start by writing the numerical flow for $\y_n +\pert$:
  \begin{equation}
   \label{pert-num-flow-1}
   \begin{array}{cl}
    \Phi_{\h}(\y_n + \pert) = \varphi_{\h} (\y_n + \pert)
    &+ \C{\op+1} \h^{\op+1} \F{\op+1}{}(\t_n,\y_n+\pert) \\
    &+ \C{\op+2} \h^{\op+2} \F{\op+2}{}(\t_n,\y_n+\pert) + O(\h^{\op+3}).
    \end{array}
  \end{equation}
We need to expand using Taylor formula of $\F{\op+1}{}$ and $\F{\op+2}{}$ in the neighborhood of $\y_n$:
  \begin{align}
  \label{Taylor-Fp}
   \F{q}{}(\t_n,\y_n + \pert) = \F{q}{}(\t_n,\y_n ) + \pert\partial_{\y} \F{q}{}(\t_n,\y_n) + \Op(\pert^2), \quad \forall \,q\geqslant \op+1,
  \end{align}
then, we use the Taylor series expansion of the exact flow around $\y_n$:
\begin{eqnarray}
 \varphi_{\h}(\y_n + \pert) &=& \y_n + \pert + \sum\limits_{j=1}^{\infty} \frac{\h^j}{j!} \f^{(j-1)}(\t_n, \y_n +\pert)\nonumber\\
 &=&   \y_n + \pert +\sum\limits_{j=1}^{\infty} \frac{\h^j}{j!} \Big( \f^{(j-1)} +\pert \partial_{\y}\big(\f^{(j-1)}\big)\Big) (\t_n, \y_n )
 + \Op(\pert^2)
 \nonumber\\
 \label{pert-exact-flow}
 &=&\varphi_{\h}(\y_n) + \pert + \pert\, \h \,\f_\y(\t_n,\y_n) + \Op(\h^{\op+3}).
\end{eqnarray}
We conclude by substituting Formulas \eqref{Taylor-Fp} and \eqref{pert-exact-flow} in \eqref{pert-num-flow-1}, then assembling terms to retrieve the desired \cref{pert-num-flow}.
\end{proof}
In the next section, we present the main Theorem stating that the imaginary part of a double composition for a numerical flow associated with a one-step method with thoroughly chosen complex coefficients is an error estimate of the approximation given by the real part.

\subsection{The main result}

For two complex constants $\gamma_1$ and $\gamma_2$, we define the double composition \cite{suzuki-1990} of the numerical flow $\Phi_\h$ associated with a one-step method, and we denote this composition by $\Psi_\h$ as follows:
\begin{equation}
\label{composition-flow}
\Psi_\h \coloneqq  \Phi_{\gamma_2\h}\circ \Phi_{\gamma_1\h}: \left\lvert
\begin{array}{ccccl}
 \mathds{R}^d &\longrightarrow &\mathds{C}^d& \longrightarrow &\mathds{C}^d\\
 \y_n &\longmapsto &\Phi_{\gamma_1\h}(\y_n) & \longmapsto & \Phi_{\gamma_2\h}\Big( \Phi_{\gamma_1\h}(\y_n) \Big)  .
 \end{array}
 \right.
\end{equation}
This means that for a given state $\y_n$ (here we choose to extend its domain of definition over the complex domain) we first apply the numerical flow for a step $\gamma_1\h$ by $\Phi_{\gamma_1\h}(\y_n)$, then use this image as input to apply again the same numerical flow for a complementary step $\gamma_2\h$ ($\gamma_1 + \gamma_2 \equiv 1$). For any complex number $z\in\mathds{C}$, we define by $\Re(z)$ its real part and $\Im(z)$ its imaginary such that:
\begin{equation*}
 z \coloneqq  \Re(z) + \i\times \Im(z),
\end{equation*}
where the complex number $\i = \sqrt{-1}$. Now we state the main result.
\begin{theorem}
 Let us take a numerical flow $\Phi_\h$ of order $\op$ and two complex coefficients $\gamma_1$ and $\gamma_2$ such that $\gamma_1 + \gamma_2=1$. We define $\Psi_\h$ as in \cref{composition-flow}. If
 \begin{equation}
\label{eq-1}
 \gamma_1^{\op+1} + \gamma_2^{\op+1}=0,
 \end{equation}
 then
 \begin{equation}
 \label{error-cond}
 \left\|\varphi_\h(\y_n) - \Re\Big(\Psi_{\h}(\y_n) \Big)\right\| \sim \Cc{}\times \left\|\Im \Big( \Psi_{\h}(\y_n)\Big) \right\|.
 \end{equation}
\end{theorem}
\begin{proof}
We start by using \cref{num-flow} to represent the first step $\Phi_{\gamma_1\h} (\y_n)$ in the composition:
\begin{eqnarray}
 \Phi_{\gamma_1\h}(\y_n) &=& \varphi_{\gamma_1\h}(\y_n) + \C{\op+1} (\gamma_1\h)^{\op+1} \F{\op+1}{}(\t_n,\y_n) + \nonumber \\
 &+& \C{\op+2} (\gamma_1\h)^{\op+2} \F{\op+2}{}(\t_n,\y_n)+ \Op(\h^{\op+3}) \nonumber\\
 \label{num-flow-1}
 &=& \varphi_{\gamma_1\h}(\y_n)  +\pert,
\end{eqnarray}
where $\pert\sim\Op(\h^{\op+1})$. Thus we can use \cref{lemma1} to write the second step in the composition as follows:
\begin{equation}
\label{num-flow-3}
\begin{array}{rl}
 \Phi_{\gamma_2\h}\Big(\Phi_{\gamma_1\h}(\y_n) \Big) &= \Phi_{\gamma_2\h}\Big( \varphi_{\gamma_1\h}(\y_n)  +\pert\Big)\\
  &=\Phi_{\gamma_2\h}\Big( \varphi_{\gamma_1\h}(\y_n)\Big) + \pert + \pert\, \gamma_2\h \,\f_\y\Big(\t_n+\gamma_1\h,\varphi_{\gamma_1\h}(\y_n)\Big)\\
 &\quad + \Op(\h^{\op+3}).
  \end{array}
\end{equation}
We finish by applying \cref{num-flow} on $\Phi_{\gamma_2\h}\Big( \varphi_{\gamma_1\h}(\y_n)\Big)$. This will produce the term $\F{\op+1}{}\big(\t_n+\gamma_1\h,\varphi_{\gamma_1\h}(\y_n)\big)$, which we recast as an expansion around $(\t_n,\y_n)$ to find the following:
\begin{equation}
\label{num-flow-2}
\begin{array}{l}
 \Phi_{\gamma_2\h}\Big( \varphi_{\gamma_1\h}(\y_n)\Big) =
 \varphi_{\gamma_2\h}\big( \varphi_{\gamma_1\h}(\y_n)\big)  \\
  \quad + \,\C{\op+1} \gamma_2^{\op+1}\h^{\op+1} \left(\left[ \F{\op+1}{} + \gamma_1\h\Big(\F{\op+1}{}\Big)^{\prime}\right](\t_n,\y_n) + \Op(\h^2) \right)  \\[6pt]
 \quad+\,  \C{\op+2} \gamma_2^{\op+2}\h^{\op+2} \left(\left[ \F{\op+2}{}  + \gamma_1\h\Big(\F{\op+2}{}\Big)^{\prime}\right](\t_n,\y_n) + \Op(\h^2) \right) \\
 \quad +\Op(\h^{\op+3}).
 \end{array}
\end{equation}
Here, functions $\big(\F{\op+1}{}\big)^{\prime}$ and $\big(\F{\op+2}{}\big)^{\prime}$ are the first order total derivatives of $\F{\op+1}{}$ and $\F{\op+2}{}$ with respect tot time $\t$. Now, we replace the Formula \eqref{num-flow-2} and the term $\pert$ presented in \eqref{num-flow-1} in relation \eqref{num-flow-3}, assemble all powers of $\h$ to find the following asymptotic expansion of the composed flow:
\begin{equation}
\begin{array}{l}
 \Phi_{\gamma_2\h}\Big(\Phi_{\gamma_1\h}(\y_n) \Big) =
 \varphi_{\gamma_2\h}\big( \varphi_{\gamma_1\h}(\y_n)\big)\\
 \quad +\h^{\op+1}\C{\op+1}\big(\gamma_1^{\op+1} + \gamma_2^{\op+1}\big)\F{\op+1}{}(\t_n,\y_n) \\[2pt]
 \quad +\h^{\op+2}\C{\op+1}\Big( \gamma_1\gamma_2^{\op+1} \big(\F{\op+1}{}\big)^{\prime}+ \gamma_2\gamma_1^{\op+1} \f_\y\cdot\F{\op+1}{}  \Big)(\t_n,\y_n)\\
 \quad +\h^{\op+2}\C{\op+2}\big(\gamma_1^{\op+2} +\gamma_2^{\op+2}\big) \,\F{\op+2}{} (\t_n,\y_n) + \Op(\h^{\op+3}).
 \end{array}
\end{equation}
We use the group property of the exact flow $\varphi_\t$ and \cref{eq-1}, to extract the real part of the composition. We mention that $\y_n$ is real, though $\varphi_\h(\y_n)$, $\big(\F{\op+1}{}\big)^{\prime}$, $\f_\y\cdot\F{\op+1}{}$ and $\F{\op+1}{}(\t_n,\y_n)$ are reals and the real part of the composed flow is written below:
\begin{equation}
\begin{array}{l}
 \Re\Big( \Psi_{\h}(\y_n)\Big) =  \varphi_\h(\y_n) \\
 \quad + \h^{\op+2} \C{\op+1}\Re \Big( \gamma_1\gamma_2^{\op+1}\Big) \big(\F{\op+1}{}\big)^{\prime}(\t_n,\y_n)\\
\quad + \h^{\op+2}\C{\op+1}\Re\Big(\gamma_2\gamma_1^{\op+1}\Big) \big(\f_\y\cdot\F{\op+1}{}\big)(\t_n,\y_n)\\
 \quad +\h^{\op+2}\C{\op+2}\Re\Big( \gamma_1^{\op+2} +\gamma_2^{\op+2} \Big)\,\F{\op+2}{} (\t_n,\y_n) + \Op(\h^{\op+3}).
\end{array}
\end{equation}
This leads us to affirm that the error between the real part and the exact flow is at least of order $\op+2$, such as:
 \begin{equation}
  \label{complex-cond}
  \varphi_\h(\y_n) - \Re\Big(\Psi_{\h}(\y_n) \Big) = \Op (\h^{\op+2}) ,
 \end{equation}
Thus,
the error is written as follows:
\begin{equation}
\label{error-comp}
\begin{array}{l}
 \e(\t_{n+1}) = \varphi_\h(\y_n) -\Re\Big( \Psi_{\h}(\y_n)\Big) = \\
 \quad - \h^{\op+2} \C{\op+1}\Re \Big( \gamma_1\gamma_2^{\op+1}\Big) \big(\F{\op+1}{}\big)^{\prime}(\t_n,\y_n)\\
\quad - \h^{\op+2}\C{\op+1}\Re\Big(\gamma_2\gamma_1^{\op+1}\Big) \big(\f_\y\cdot\F{\op+1}{}\big)(\t_n,\y_n)\\
 \quad -\h^{\op+2}\C{\op+2}\Re\Big( \gamma_1^{\op+2} +\gamma_2^{\op+2} \Big)\,\F{\op+2}{} (\t_n,\y_n) + \Op(\h^{\op+3}).
\end{array}
\end{equation}
For the asymptotic error \eqref{error-comp}, we have three cases where each case corresponds to which term among $\big(\F{\op+1}{}\big)^{\prime}$, $\f_\y\cdot\F{\op+1}{}$ or $\F{\op+2}{}$ is the leading one. Before continuing, we need the following identity. It proof is proven in \ref{app1}.
\begin{equation}
\label{Im_12}
\Im\big(\gamma_1^{\op+2} +\gamma_2^{\op+2}\big)=0.
\end{equation}
Thus, we write the imaginary part of the composed flow as follows:
\begin{equation}
\label{imag-comp}
 \begin{array}{l}
  \Im\left(\Psi_{\h}(\y_n) \right) = \h^{\op+2} \C{\op+1}\Im \Big( \gamma_1\gamma_2^{\op+1}\Big) \big(\F{\op+1}{}\big)^{\prime}(\t_n,\y_n)\\
  \quad + \h^{\op+2}\C{\op+1}\Im\Big(\gamma_2\gamma_1^{\op+1}\Big) \big(\f_\y\cdot\F{\op+1}{}\big)(\t_n,\y_n)+ \Op(\h^{\op+3}).
 \end{array}
\end{equation}
If the leading term of the error in \eqref{error-comp} is $\big(\f_\y\cdot\F{\op+1}{}\big)(\t_n,\y_n)$, then, the norm of the imaginary part \eqref{imag-comp} associated to $\Psi_\h$ is asymptotically equivalent to the norm of the error \eqref{error-comp} by a constant given by:
\begin{equation*}
\Cc{1} \coloneqq {\Re\Big(\gamma_2\gamma_1^{\op+1}\Big)/}{\Im\Big(\gamma_2\gamma_1^{\op+1}\Big)}.
\end{equation*}
If $\F{\op+2}{}$ is the leading term in \eqref{error-comp},
and if $\big(\F{\op+1}{}\big)^{\prime}\sim \F{\op+2}{}$, which is true for a certain number of schemes, then the constant would be
\begin{equation*}
\Cc{2}\ \coloneqq \ {\Re\Big(\gamma_1\gamma_2^{\op+1}\Big)}/{\Im\Big(\gamma_1\gamma_2^{\op+1}\Big)} + {\C{\op+2}\Re\Big( \gamma_1^{\op+2} +\gamma_2^{\op+2} \Big)}/{\C{\op+1}\Im \Big( \gamma_1\gamma_2^{\op+1}\Big)}.\end{equation*}
Since in practice we do not have any more accurate idea about $\C{\op+1}$ and $\C{\op+2}$, we consider that their ratio $\C{\op+2}/\C{\op+1} \approx \op+2$.
If $\big(\F{\op+1}{}\big)^{\prime}$ is the leading term, then the constant will be also equal to $\Cc{2}$. To this end, we take the constant $\Cc{}$ to be
\begin{equation}
\label{max_C1_C2}
\Cc{} = \max(\Cc{1},\Cc{2}).
\end{equation}
\end{proof}

\subsection{Discussions}
\cref{complex-cond} states that the real part of the composition technique approximates the solution up to order $\op+1$. \cref{error-cond} states that the imaginary part of the output of the composition is an error estimate of the real part of $\Psi_\h$ of order $\op+1$. One has a numerical flow associated to a one-step method of order $\op$, this is of importance when dealing with \ac{ATS}. The fact of composing it twice with complex coefficients will not only give us a higher-order numerical approximation in the real part, but it will produce also an error estimate of order $\op+1$ in the imaginary part. Therefore, this error estimate is effectively used in the adaptive time step to follow up the dynamics of the numerical simulation and update the time step according to a user tolerance. For safety measures, we update the time step as follows:
\begin{equation}
 \label{update-ts}
 \h_{n+1} \coloneqq C \times  \,\h_n\left(\frac{\tol}{\Cc{}\|\Im\left(\Psi_{\h}(\y_n) \right)\|}\right)^{\frac{1}{\op+1}},
\end{equation}
with  $C \approx 0.9$
and the $\left\| \cdot \right\|$ is considered to be the Euclidean norm.

\subsection{Algorithm}
\label{sec-alg}
Hereafter, we present the steps to follow in order to produce a numerical simulation, starting with initial condition $\y_0$, with adaptive time step using a one-step numerical flow $\Phi_\h$ of order $\op$.

\begin{algorithm}[ht!]
 \caption{Adaptive numerical simulation using $\Phi_{\h}$\label{alg-adap}}
\begin{algorithmic}
\Require $\f$, $\t_0$, $\y_0$, $\T{}$, $\h_0$, $\tol$, $\Phi_{\h}$, $\op$
\State
\begin{flalign}
\label{gamma_i}
\gamma_1 &\gets \frac{1}{2} + \frac{\i}{2}\frac{\sin(\frac{\pi}{\op+1})}{1+\cos(\frac{\pi}{\op+1})} &&
\end{flalign}
\State $\gamma_2 \gets 1-\gamma_1$
\State $\Cc\, \gets \max\big(\Cc{1},\Cc{2}\big)$
\State $\h_n \gets \h_0$
\State $\t_n \gets \t_{0} +\h_n$
\State $\y_n \gets \y_0$
\While{$\t_n\leqslant \T{}$}
\State $\y_{n+1,1} \gets \Phi_{\gamma_1 \h_n}(\y_n)$
\State $\y_{n+1,2} \gets \Phi_{\gamma_2 \h_n}(\y_{n+1,1})$
\State $\y_{n+1} \gets \Re\big(\y_{n+1,2}\big)$
\State $\e_{n+1} \gets \Cc{}\times \left\|\Im\left( \y_{n+1,2} \right) \right\|$
\State \begin{flalign}
        \label{adapt-form}
        \quad \h_{n+1} &\gets C\times \displaystyle\h_n\left( \frac{\tol}{\e_{n+1}}\right)^{\frac{1}{\op+1}} &&
       \end{flalign}
\State $\t_{n+1} \gets \t_n + \h_{n+1}$
\State $n\gets n+1$
\EndWhile
\end{algorithmic}
\end{algorithm}

Regarding the computational cost, this technique generate additional costs due to complex arithmetic. It will be shown in the next section that, despite of these additional costs, they are negligible comparing the increasing of the order of accuracy and having error estimate in the imaginary part. We will show that, composed schemes produce numerical approximation with higher accuracy and lower computational costs regarding the basic ones.

\section{Composition of familiar numerical schemes: Linear stability and rate of convergence}
\label{sec3}

We consider here some numerical integrators and their double composition. A comparison between the linear stability domain of the basic integrator and the composition will be shown. To do that, we start with the linear equation $\y^\prime = \lambda \y$, and check the domain $\D = \{z\in \mathbb{C},  \,\text{s.t.} \, |\P(z)|\leqslant 1 \}$, where the function $\P(z)$ is defined from a given numerical integrator $\Phi_{\h}(\y) = \P(z)\,\y$, where $z = \h\lambda$. We denote by $\P(z)$ the stability polynomial of the associated numerical integrator $\Phi_\h$. We compare this domain with the stability domain of the resulted composed scheme $\Re\big( \Phi_{\gamma_2\h}\circ\Phi_{\gamma_1\h} \big)$, denoted here by $\Dc$ and resulted form the associated stability polynomial $\Pc(z)$ as follows:
\begin{equation}
\Pc(z) \equiv  \Re\big( \P(\gamma_2 z)\cdot\P(\gamma_1 z)\big).
\end{equation}
Therefore the linear stability domain is the following part in the complex plane:
\begin{equation*}
 \Dc = \{z\in \mathbb{C},  \,\text{s.t.} \, |\Pc(z)|\leqslant 1 \}.
\end{equation*}
For a set of values $N_i$, we define a set of time steps $\h_i = \frac{\T{}}{N_i}$, every scheme and its composition will be performed to solve the following equation $\frac{\d\y}{\d\t} = -\y^3$ with the initial condition $\y(0)=1$ over the interval $[0,\T{}]$ with $\T{}=2$. The time interval will be decomposed uniformly on $N_i$ points and solution will be approximated over a discrete set of instants $\t_n = n\h$ for $n \in \{1,\ldots,N_i\}$. The error at every instant $\t_n$ is calculated with the exact solution $\y(\t) = \frac{1}{\sqrt{1+2\cdot \t}}$. The global error is approximated by $\overline{\e}_{\h_i}$ using the trapezoidal method:
\begin{equation}
 \overline{\e}_{\h_i} \coloneqq \sum\limits_{n=1}^{N_i-1} \h_i\cdot \left\vert \y_n - \y(\t_n) \right\vert + \frac{\h_i}{2}\cdot \left\vert \y_{N_i} - \y(\T{}) \right\vert
\end{equation}
Now, we define the rate of convergence as follows:
\begin{equation}
 \label{ROC}
 \mathrm{ROC} \coloneqq \lim\limits_{\substack{\h_i,\h_j \rightarrow 0\\\h_j>\h_i}} \frac{\log_{10}\Big(\displaystyle\frac{\overline{\e}_{\h_i}}{\overline{\e}_{\h_j}} \Big)}{\log_{10}\Big(\displaystyle\frac{{\h_i}}{{\h_j}} \Big)}.
\end{equation}
which is evaluated for the set of values $\h_i$ and results are presented in \cref{tab5}. Values in every row should converge to $\op$ when using a basic integrator $\Phi_{\h}$ of order $\op$, and to a value $\geqslant\op+1$ when $\Re \big(\Psi{\h}\big)$ is used. We confirm the fact that composition technique presented in this paper increases the order of the integration as presented in the following figures and \cref{tab5}.

\subsection{\ac{RK} families}
\label{sec312}
The family of \ac{RK} schemes is the most widely used schemes in numerical simulation due to their theoretical foundation and their versatile applications in many fields. They were initially proposed by Runge \cite{Runge_1895} and are classified by the number of stages $s$. The general form of a \ac{RK} scheme is given as follows:
\begin{align*}
	 \y_{n+1} &\coloneqq \Phi^{\RK{s}}_{\h_n}(\y_n)  \equiv \y_n + \h_n\sum\limits_{i=1}^{s} b_i K_i, \\
	K_i &\coloneqq \f\Big( t_n + c_i \h_n , \y_n + \sum\limits_{j=1}^{s}a_{ij} K_j\Big), \\
    c_i &\coloneqq \sum\limits_{j=0}^sa_{ij}.
\end{align*}
Thanks to the works of Butcher \cite{butcher-1969, book:butcher}, the \ac{RK} schemes are represented by the following table called the Butcher tableau:
\[
\renewcommand\arraystretch{1.2}
\begin{array}
	{c|ccc}
	c_1 &a_{11}  & \ldots &a_{1s}\\
	\vdots &  \vdots & & \vdots \\
	c_s&a_{s1}& \ldots &a_{ss} \\
	\hline
	& b_1 & \ldots &b_s
\end{array}
\qquad\qquad
\begin{array}{c|c}
 \mathbf{c} & \mathbf{A}\\\hline
 & \mathbf{b}
\end{array}
\]
Implicit and explicit \ac{RK} schemes are considered to study their linear stability after composition. According to \cite[page 44]{book:hairer2}, the stability polynomial of a \ac{RK} scheme is defined by the following:
\begin{equation}
 \label{poly_Stab_RK}
 \P(z) \coloneqq 1 + z\cdot \Big(\mathbf{b}^\top\cdot \big(\I - z\cdot\mathbf{A}\big)^{-1}\cdot E\Big),
\end{equation}
where $E$ is the matrix of ones and $\I$ is the identity matrix.
To this end, \ac{ERK} methods are defined when a second set of coefficients $\mathbf{b}^*$ is added to a given \ac{RK} scheme of order $\op$ defined by its Butcher Tableau ($\mathbf{c,A,b}$). This allows to evaluate a second approximation $\hy_{n+1} = \sum\limits_{i=1}^s{b_i^* K_i}$ of order $\op-1$. Thus an error estimate could be evaluated as $\e^*_n \,\coloneqq\ \y_{n+1} - \hy_{n+1} $. 

\subsubsection{First illustration: Composition of the first order forward Euler scheme}
The basic explicit integrator of order one is the Euler scheme given by the numerical flow $\Phi^{\FE}_{\h_n}$ such that the image of a given $\y_n$ is presented below:
\begin{equation*}
	\Phi^{\FE}_{\h_n}(\y_n) \ \coloneqq \ \y_n + \h_n \f(\t_n, \y_n).
\end{equation*}
Its double composition prescribes the following coefficients: $ \gamma_1 = 1/2 + \i/2 $ and $\gamma_2 = 1/2 - \i/2 $, where \cref{fig:illustration_scheme_s_3} illustrates it. In the first step, we denote by ${\y_{n+1,1}} $ the image of the flow with the step $\gamma_1\h_n$ and is computed as follows:
\begin{align*}
  \y_{n+1,1}\ \coloneqq \ \Phi^{\FE}_{\gamma_1\h_n}(\y_n) &= \y_n + \gamma_1 \h_n\f\big(\t_n,\y_n\big) =\y_n + \frac{\h_n}{2}\f\big(\t_n,\y_n\big) +  \i\frac{\h_n}{2} \f\big(\t_n,\y_n\big)\,.
\end{align*}
Next, we evaluate $\f(\t_n+\gamma_1\h_n,\y_{n+1,1}) \equiv X+\i Y$, which is an approximation of  $\y^{\prime}(t_{n}+\gamma_1\h_n)$, and multiply it by $\gamma_2\h_n$ in order to compute the second step $\Phi^{\FE}_{\gamma_2\h_n}(\y_{n+1,1})$ denoted by $\y_{n+1,2}$:
\begin{figure}[h]
	\begin{center}
		\includegraphics[height=6.5cm]{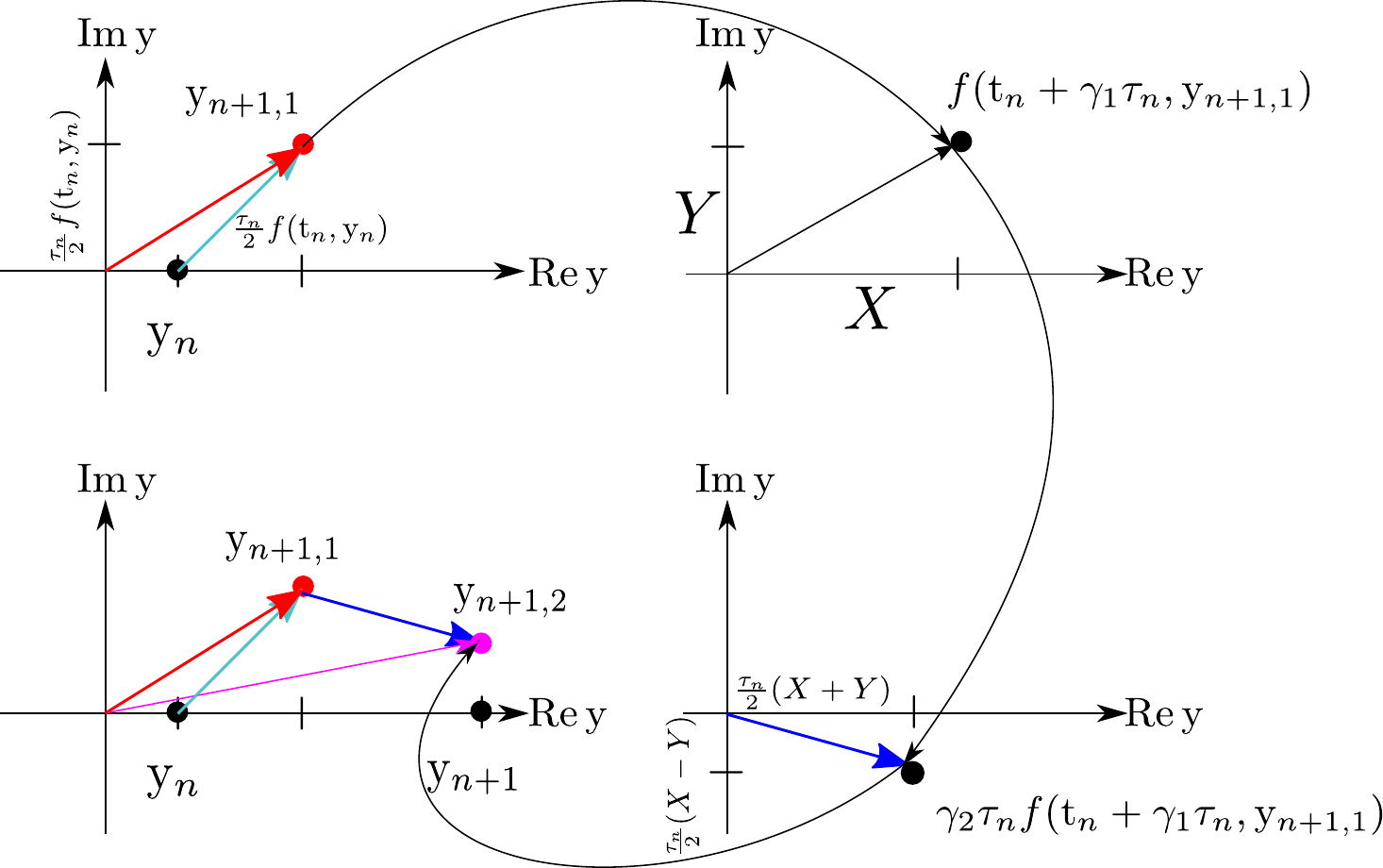}
	\end{center}
	\caption{Illustration of composing twice the Forward Euler scheme.}
	\label{fig:illustration_scheme_s_3}
\end{figure}
\begin{align*}
	 \y_{n+1,2} \ \coloneqq \ & \Phi^{\FE}_{\gamma_2\h_n}(\y_{n+1,1}) \equiv \y_{n+1,1} + \gamma_2 \h_n\f\big(\t_n+\gamma_1 \h_n, \y_{n+1,1}\big)\\
	= &\y_{n+1,1} + \frac{\h_n}{2}\f\big(\t_n+\gamma_1 \h_n,\y_{n+1,1}\big) - \i\frac{\h_n}{2}\f\big( \t_n+\gamma_1 \h_n,\y_{n+1,1}\big)\,.
\end{align*}
By replacing $\y_{n+1,1}$, we end by presenting the real and imaginary part of the composition:
\begin{align*}
	\Phi^{\FE}_{\gamma_2\h_n} \circ \Phi^{\FE}_{\gamma_1\h_n} (\y_{n}) &= \Big(\y_n  + \frac{\h_n}{2}\Big[\f\big(\t_n,\y_n\big) + X +Y\Big]\Big) + \i\,\frac{\h_n}{2} \Big( \f\big(\t_n,\y_n\big) +Y -X  \Big)\\
	& = \Re \big(\y_{n+1,2}\big) + \i\,\Im\big({\y_{n+1,2}}\big)\,,
\end{align*}
where the real part is the approximation of the solution at $\t=\t_n+\h_n$, and the imaginary part will be used as an error estimator to adapt locally the time step. We denote by $\Phi_\h^{\CRK1}(\y_n) \coloneqq \Re\big( \Phi^{\FE}_{\gamma_2\h_n} \circ \Phi^{\FE}_{\gamma_1\h_n} (\y_n)\big)$ the \ac{CRK} scheme of first order.

\begin{figure}[!ht]
\centering
 \begin{subfigure}{0.3\textwidth}
 \includegraphics[width=\textwidth]{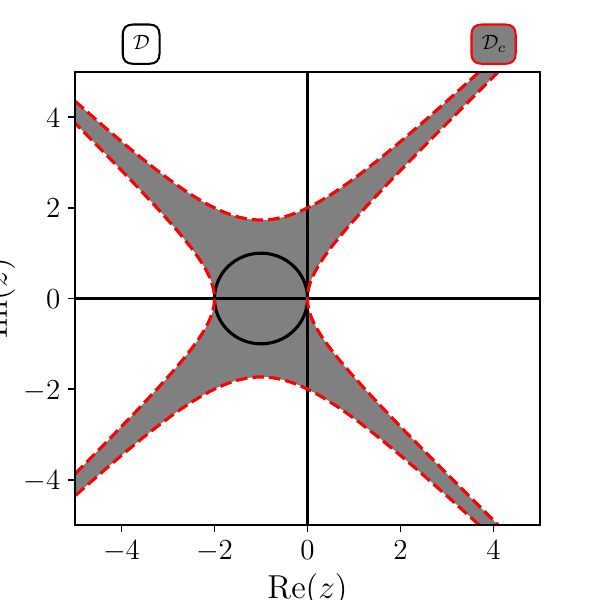}
 \caption{$\D$ (inside black line) and $\Dc$ (the grey region). \label{FE1}}
 \end{subfigure}
 \begin{subfigure}{0.3\textwidth}
 \includegraphics[width=\textwidth]{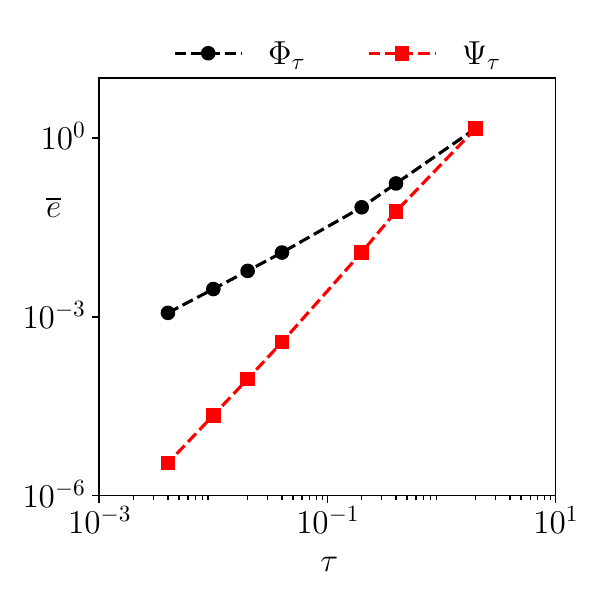}
 \caption{Global error $\overline{\e}_\h$ versus time step $\h$. \label{FE2}}
 \end{subfigure}
 \begin{subfigure}{0.3\textwidth}
 \includegraphics[width=\textwidth]{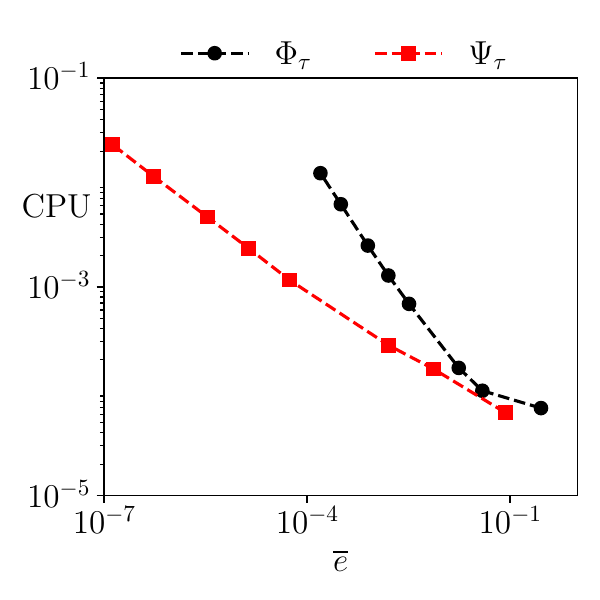}
 \caption{CPU versus global error $\overline{\e}$.\label{FE3}}
 \end{subfigure}
 \caption{Comparison between the Linear stability region of the numerical flow $\Phi_\h^{\FE}$ of the Forward Euler scheme and its double composition $\Phi_\h^{\CRK1}$ (left panel), the order of convergence (middle) and their CPU time (right panel). \label{FE}}
 \end{figure}

In \cref{FE1}, we present the stability region $\D$ of the basic integrator, surrounded by the black solid line, and $\Dc$ of its composition represented by the region coloured in grey and surrounded by the dashed red line. It is clear that $\D \subset \Dc$. However, both domains have the same boundary point crossing the negative real axis.
In \cref{FE2} the rate of convergence for both, the basic integrator and its composition,  are plotted. It is shown that the rate of convergence is improved by the composition, as same as the computational time (see \cref{FE3}). For additional details, the CPU time needed to achieve the simulation with global precision $\overline{\e}$ of order $10^{-3}$ by the composed flow $\Phi_\h^{\CRK1}$ is six time smaller (see \cref{tab4}) then using the basic integrator $\Phi_\h^{\RK1}$. This ratio increases with higher prescription of precision.

\subsubsection{Second-order scheme}
We present here the Butcher tableau of the second-order explicit \ac{RK} scheme:
\begin{table}[!ht]
 \centering
  \caption{Butcher tableau of \ac{RK}2 method. \label{tab1}}
    \renewcommand\arraystretch{1.2}
\begin{tabular}{c|cc}
	0\\
	$\alpha$ & $\alpha$\\
	\hline
	& $1-\frac{1}{2\alpha}$ & $\frac{1}{2\alpha}$
\end{tabular}
\end{table}
For $\alpha = 1/2$, we have the midpoint method, and the Heun's method is recovered for $\alpha = 1$. In both cases, the coefficients of composition are: $\gamma_{1} = 1/2 + \i\sqrt{3}/6$ and $\gamma_{2} = 1/2 - \i\sqrt{3}/6$. We denote by \ac{CRK}2 scheme the composition of the second order \ac{RK}2.

\begin{figure}[!ht]
\centering
 \begin{subfigure}{0.3\textwidth}
 \includegraphics[width=\textwidth]{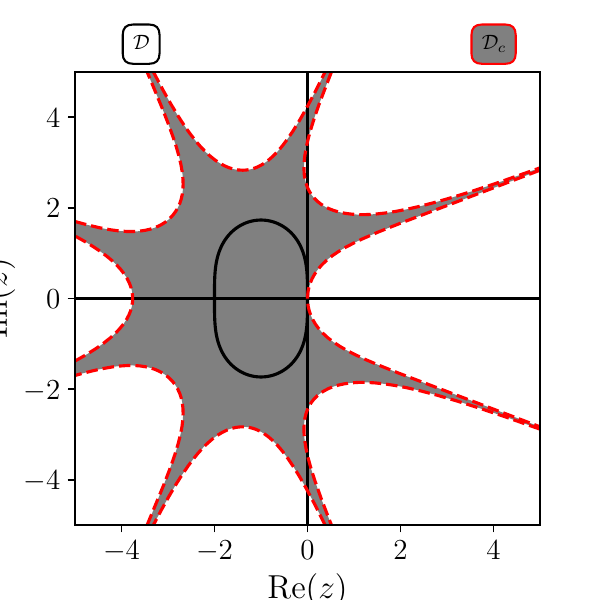}
 \caption{$\D$ (inside black line) and $\Dc$ (the grey region). \label{RK21}}
 \end{subfigure}
 \begin{subfigure}{0.3\textwidth}
 \includegraphics[width=\textwidth]{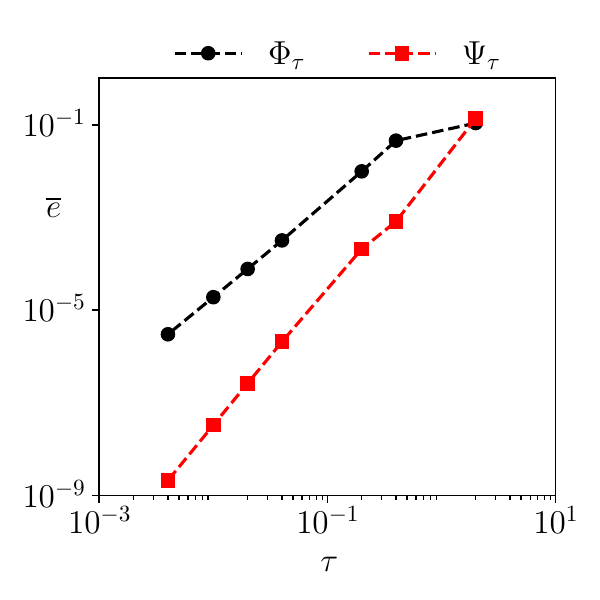}
 \caption{Global error $\overline{\e}$ versus time step $\h$. \label{RK22}}
 \end{subfigure}
 \begin{subfigure}{0.3\textwidth}
 \includegraphics[width=\textwidth]{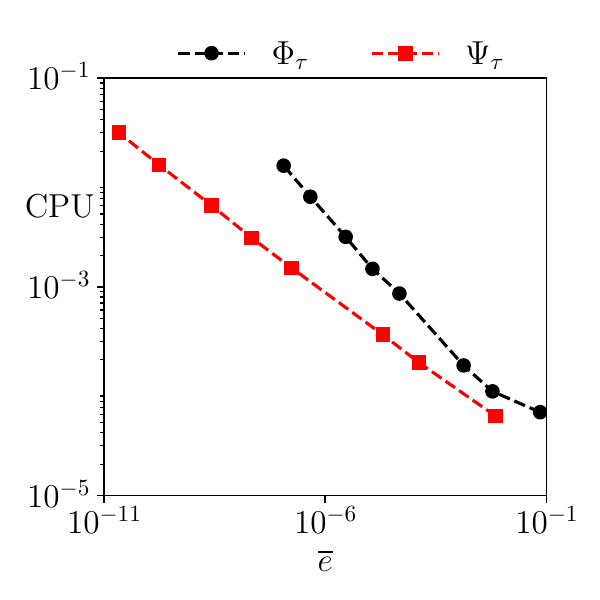}
 \caption{CPU versus global error $\overline{\e}$.\label{RK23}}
 \end{subfigure}
 \caption{Comparison between the Linear stability region of the numerical flow of the \ac{RK}2 scheme with $\alpha=1/2$ and its double composition (left panel), the order of convergence (middle) and their CPU time (right panel). \label{RK2}}
 \end{figure}

In this scheme, we see also an increasing in the linear stability domain as the boundary point of $\Dc$ that is crossing the negative real axis has an absolute value twice bigger than the one of $\D$. The order of convergence is also increased by one (see \cref{RK22} and \cref{tab5}). The computational efficiency is improved by the composition, where the CPU time needed to achieve the simulation with global precision $\overline{\e}$ of order $10^{-5}$ by the composed flow $\Re\big(\Psi_\h^{\CRK2}\big)$ is five time smaller (see \cref{tab4}) then using the basic integrator $\Phi_\h^{\RK2}$. \cref{RK23} present the CPU for a range of global precisions.

\subsubsection{Fourth-order classical \ac{RK}}
This is a scheme composed of four stages and it is also a fourth-order scheme.
Coefficients of composition are evaluated using the formula \eqref{gamma_i} for $\op=4$:
In this case, the coefficients of the composition are given approximately by: $\gamma_1 \coloneqq \frac{1}{2} + \i\,\frac{\sin(\pi/5)}{1 + \cos(\pi/5)} \approx 0.5 + \i\times0.3249196964 $ and $\gamma_2 \coloneqq \frac{1}{2} - \i\,\frac{\sin(\pi/5)}{1 + \cos(\pi/5)}$.
\begin{figure}[!ht]
\centering
 \begin{subfigure}{0.3\textwidth}
 \includegraphics[width=\textwidth]{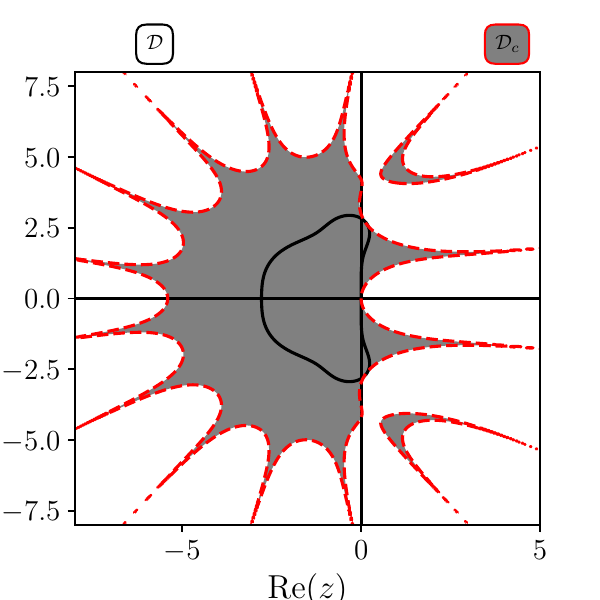}
 \caption{$\D$ (inside black line) and $\Dc$ (the grey region). \label{RK41}}
 \end{subfigure}
 \begin{subfigure}{0.3\textwidth}
 \includegraphics[width=\textwidth]{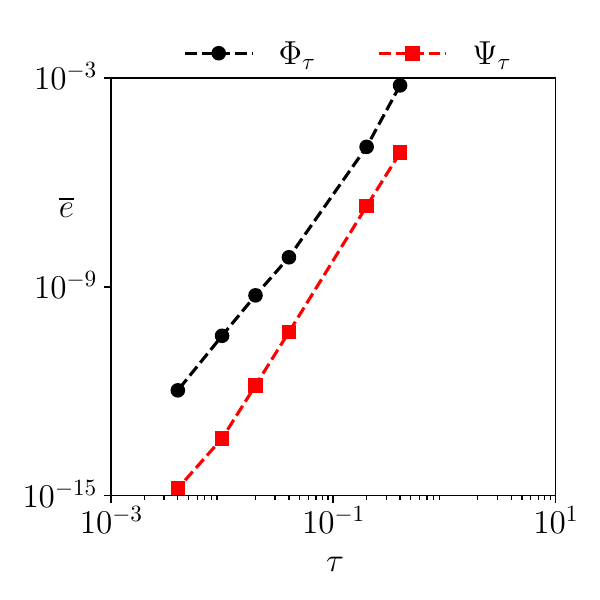}
 \caption{Global error $\overline{\e}$ versus time step $\h$. \label{RK42}}
 \end{subfigure}
 \begin{subfigure}{0.3\textwidth}
 \includegraphics[width=\textwidth]{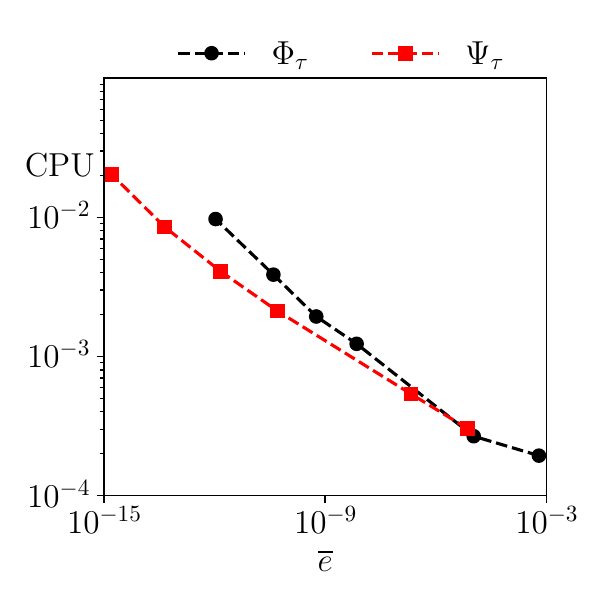}
 \caption{CPU versus global error $\overline{\e}$.\label{RK43}}
 \end{subfigure}
 \caption{Comparison between the Linear stability region of the numerical flow of the \ac{RK}4 scheme and its double composition (left panel), the order of convergence (middle) and their CPU time (right panel). \label{RK4}}
 \end{figure}
 The comparison between the basic integrator $\Phi_\h^{\RK4}$ and its double composition $\Phi_\h^{\CRK4}$ is presented in \cref{RK4}. The left panel shows an increase in the linear stability after composition as seen by the boundary points of $\D$ and $\Dc$ crossing the negative real axis.
 The middle one shows an improvement of the precision when the simulation is done with the same value of the time step. The curves in this figure could not show graphically how the composition increase the rate of convergence.  This is illustrated by the rate of convergence presented in \cref{tab5}. The right panel compares the performance between the $\Phi_\h$ and $\Psi_\h$, where in \cref{tab4}, the composition reduces by a factor $1.46$ the computational time if the target accuracy is of order $10^{-9}$.

\subsubsection{Fourth-order \ac{GRK} scheme}
\label{sec-GRK2}
\begin{figure}[!ht]
\centering
 \begin{subfigure}{0.3\textwidth}
 \includegraphics[width=\textwidth]{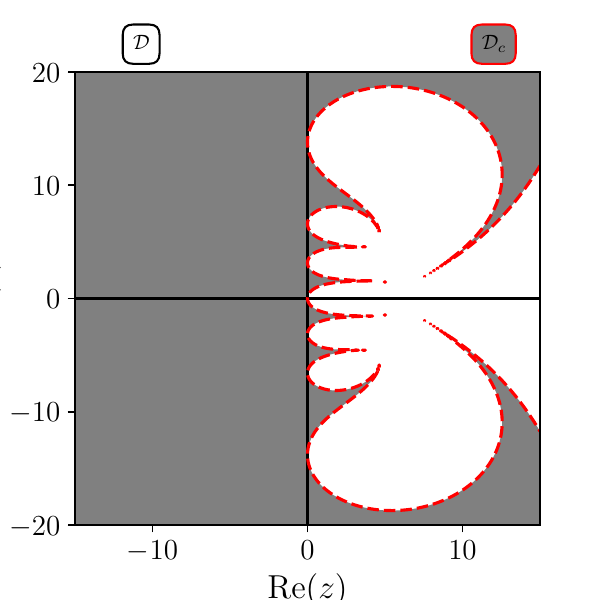}
 \caption{$\D$ (inside black line) and $\Dc$ (the grey region). \label{GRK21}}
 \end{subfigure}
 \begin{subfigure}{0.3\textwidth}
 \includegraphics[width=\textwidth]{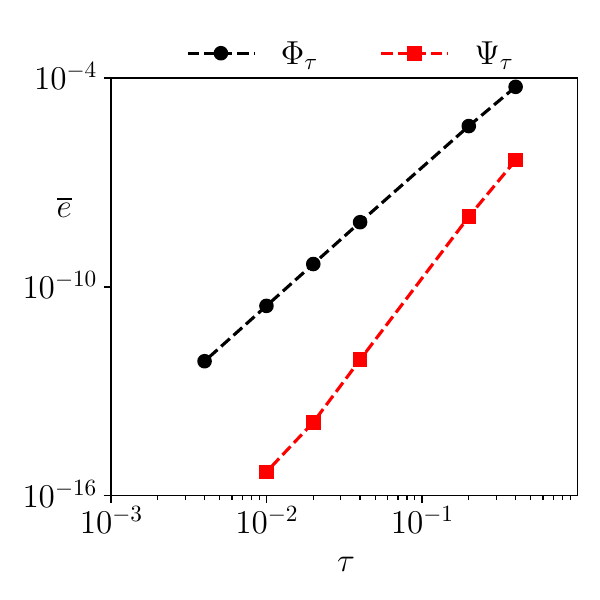}
 \caption{Global error $\overline{\e}$ versus time step $\h$. \label{GRK22}}
 \end{subfigure}
 \begin{subfigure}{0.3\textwidth}
 \includegraphics[width=\textwidth]{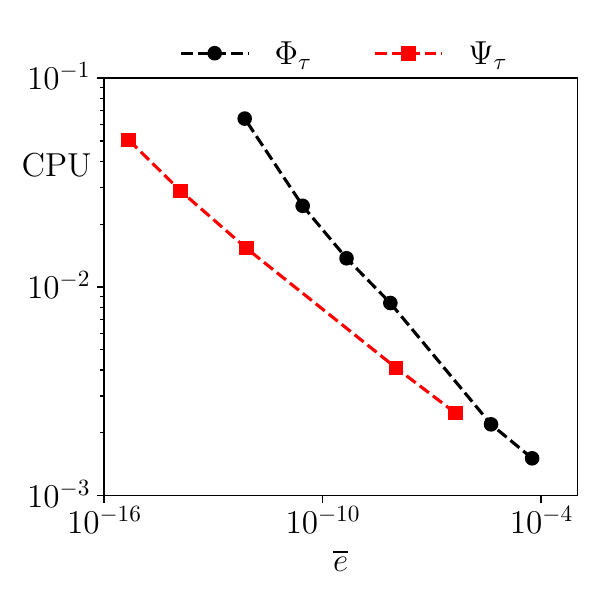}
 \caption{CPU versus global error $\overline{\e}$.\label{GRK23}}
 \end{subfigure}
 \caption{Comparison between the Linear stability region of the numerical flow of the \ac{GRK}2 scheme and its double composition (left panel), the order of convergence (middle) and their CPU time (right panel). \label{GRK2}}
 \end{figure}

\ac{GRK} methods are collocation methods. Any collocation method is based on finding a polynomial $p(\t)$ of order $s$ verifying the following relations:
\begin{equation*}
	p(\t_n) = \y_0, \quad p^{\prime}(\t_n +c_i \h) = \f\big(\t_n +c_i \h, p(\t_n +c_i \h)\big),\quad i\in \S{1}^{s},
\end{equation*}
where $\{c_i\,\vert\, i \in\S{1}^{s}\}$ are distinct real numbers (usually chosen in $[0,1]$ ). In this way, the solution is approximated by $\y_{n+1} \coloneqq p(\t_n+ \h)$. It was proven that the collocation methods are equivalent to $s$-stages \ac{RK} method, where coefficients $a_{ij}$ and $b_i$ are evaluated using the Lagrange polynomials $\l_j$ as follows:
\begin{align*}
	l_j(\h)& \ \coloneqq \ \prod\limits_{\substack{i=1\\i\neq j}}^{s}\frac{(\h-c_i)}{(c_j-c_i)}, & a_{ij} &\coloneqq \int_0^{c_i} l_j(x) \d x,  & b_i &\coloneqq \int_0^1 l_i(x) \d x.
\end{align*}
Following the above definition, if the collocation points $\{c_i\,\vert\, i \in\S{1}^{s}\}$ are taken to be points of the $s$\up{th} shifted Gau\ss-Legendre polynomials, we find the \ac{GRK} method with $s$ stages are of order $2s$ \cite[Section II.1.3]{hairer2002geometric}. In the case of considering $s=2$, the \ac{GRK} has two stages and is defined by Butcher tableau as follows:
\begin{table}[!ht]
 \centering
  \caption{Butcher tableau of \ac{GRK}2 method. \label{tab2}}
    \renewcommand\arraystretch{1.2}
\begin{tabular}{c|cc}
	 $\frac{1}{2}-\frac{\sqrt3}{6}$ &  $\frac{1}{4}$ & $\frac{1}{4}-\frac{\sqrt3}{6}$  \\
	$\frac{1}{2}+\frac{\sqrt3}{6}$  & $\frac{1}{4}+\frac{\sqrt3}{6}$ & $\frac{1}{4}$ \\
	\hline
	& $\frac{1}{2}$ & $\frac{1}{2}$
\end{tabular}
 \end{table}
This is a fourth-order numerical scheme, for which coefficients of composition are the same as given above in the part of the fourth-order classical \ac{RK} scheme. We present in \cref{GRK2} the comparison between the scheme and its double composition. The \ac{GRK}2 is A-stable as the domain of stability $\D$ the left part of the complex plane (see \cref{GRK21}). The composition is also A-stable as the domain of stability $\Dc$ associated to $\Phi^{\CGRK2}_\h$ contains $\D$. \cref{GRK22} plots the global error $\overline{\e}$ versus different time steps $\h$, where it is clear that the composition increase by one the order. \cref{GRK23} shows that the composition helps reducing the time of computation as the CPU needed to achieve the simulation for a given precision is smaller than the one done by the basic integrator.

\subsubsection{Lobatto IIIA method with three stages}

\begin{figure}[!ht]
\centering
 \begin{subfigure}{0.3\textwidth}
 \includegraphics[width=\textwidth]{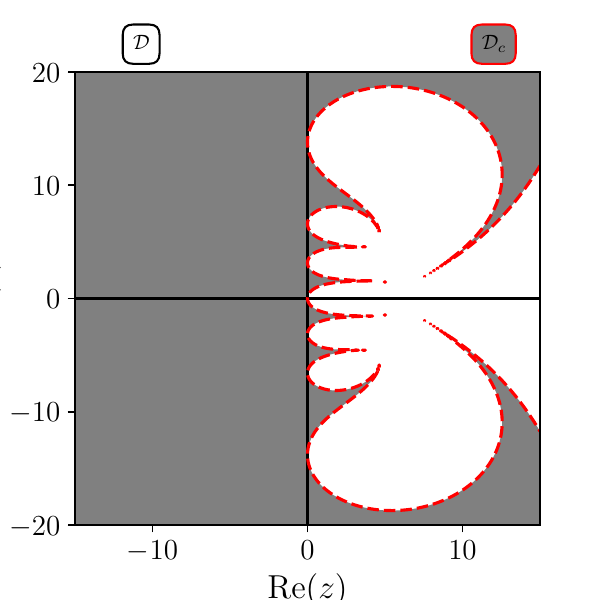}
 \caption{$\D$ (inside black line) and $\Dc$ (the grey region). \label{LIIIA31}}
 \end{subfigure}
 \begin{subfigure}{0.3\textwidth}
 \includegraphics[width=\textwidth]{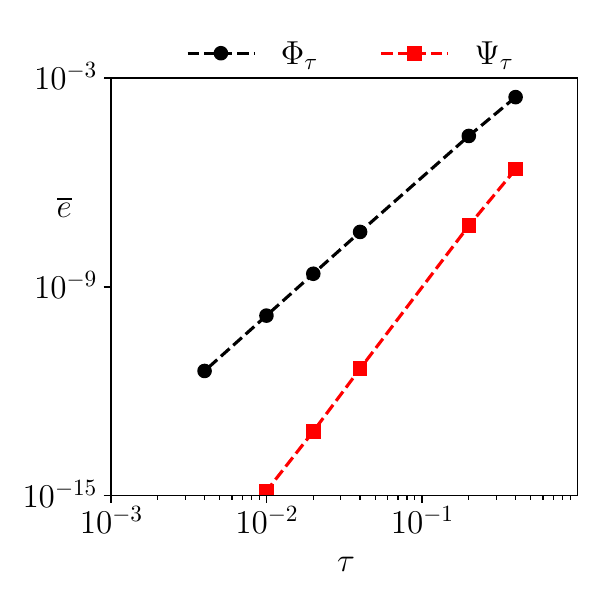}
 \caption{Global error $\overline{\e}$ versus time step $\h$. \label{LIIIA32}}
 \end{subfigure}
 \begin{subfigure}{0.3\textwidth}
 \includegraphics[width=\textwidth]{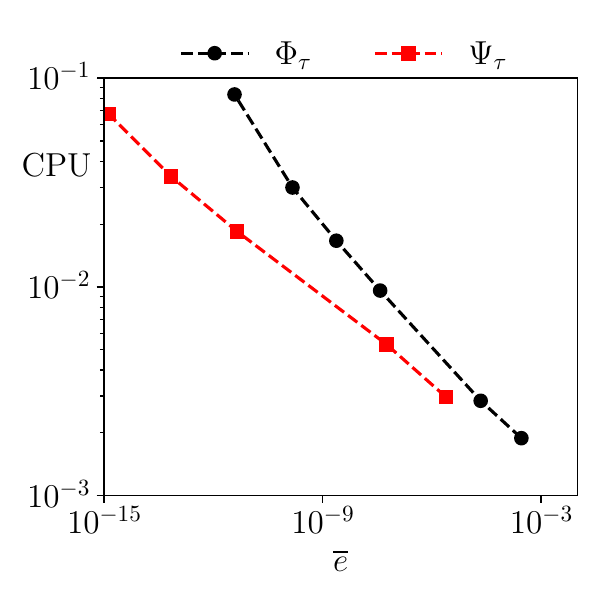}
 \caption{CPU versus global error $\overline{\e}$.\label{LIIIA33}}
 \end{subfigure}
 \caption{Comparison between the Linear stability region of the numerical flow of the Lobatto IIIA scheme of order four and its double composition (left panel), the order of convergence (middle) and their CPU time (right panel). \label{LIIIA3}}
 \end{figure}

 We end this section by showing the comparison in \cref{LIIIA3} relative to the Lobatto IIIA method given by its Butcher \cref{tab3}.  It is a scheme with three stages but is of order four.

 \begin{table}[!ht]
 \centering
  \caption{Butcher tableau of Lobatto IIIA method. \label{tab3}}
    \renewcommand\arraystretch{1.2}
\begin{tabular}{c|ccc}
   $0$ & $0$ & $0$ & $0$ \\
   $\frac{1}{2}$ & $\frac{5}{24}$ & $\frac{1}{3}$ & -$\frac{1}{24}$\\
   $1$ & $\frac{1}{6}$& $\frac{2}{3}$& $\frac{1}{6}$\\
   \hline
   &  $\frac{1}{6}$& $\frac{2}{3}$& $\frac{1}{6}$
  \end{tabular}
 \end{table}

 To stress on the CPU improvement, \cref{tab4} presents the ratio, for all the above schemes,  of the CPU of $\Phi_\h$ with the CPU of its double composition given by $\Re\big( \Psi_\h\big)$. In the case of the composition of the Lobatto IIIA3, the time computational can be reduced by a factor of 2.32 for a prescribed precision of order $10^{-9}$.

 \begin{table}[!ht]
 \centering\renewcommand\arraystretch{1.2}
 \caption{$\mathrm{ROC}$ for basic integrators and their double compositions.\label{tab5}}
\begin{tabular}{cc|cccc}
&$\h$ &   0.2& 0.04& 0.02 & 0.01  \\\hline
\multirow{2}{*}{\ac{RK}1}& $\Phi_\h$ & 1.328&1.086&1.023&1.011\\
 & $\Re(\Psi_\h)$& 2.277&2.156&2.051&2.026\\\hline
\multirow{2}{*}{\ac{RK}2}& $\Phi_\h$ & 2.200&2.132&2.045&2.023 \\
 & $\Re(\Psi_\h)$& 1.977&2.859&2.984&2.994\\\hline
\multirow{2}{*}{\ac{RK}4}& $\Phi_\h$ & 5.873&4.544&3.639&3.861\\
 & $\Re(\Psi_\h)$& 5.105&5.182&5.092&5.068\\\hline
\multirow{2}{*}{GR2}& $\Phi_\h$ & 3.748&3.955&3.996&3.999\\
 & $\Re(\Psi_\h)$& 5.402&5.882&5.983&5.839\\ \hline
\multirow{2}{*}{Lobatto IIIA 4}& $\Phi_\h$ & 3.710&3.949&3.996&3.999\\
& $\Re(\Psi_\h)$& 5.404&5.883&5.989&5.996\\\hline
\end{tabular}
 \end{table}

 We conclude that the higher the prescribed ratio is, the better the use of the composition technique is in reducing time computation.
 We add also the ROC to every basic integrator and its double composition in \cref{tab5}.

 \begin{table}[!ht]
 \centering\renewcommand\arraystretch{1.2}
 \caption{Ratio of the CPU needed by the basic flow $\Phi_\h$ over the CPU needed by its double composition $\Re\big(\Psi_\h\big)$ to perform a simulation with a global precision $\overline{\e}$ for the above defined schemes.\label{tab4}}
\begin{tabular}{c|cccccc}
$\overline{\e}$ & $10^{-3}$ & $10^{-5}$ & $10^{-7}$ & $10^{-9}$ & $10^{-12}$ & $10^{-15}$\\ \hline
\ac{RK}1&6.96&33.80&164.09&796.71&8523.74&91193.01\\\hline
\ac{RK}2&3.00&4.54&6.87&10.41&19.39&36.13\\\hline
 \ac{RK}4&1.24&1.31&1.39&1.46&1.59&1.72\\\hline
\ac{GRK}2&1.28&1.52&1.80&2.14&2.77&3.58\\ \hline
 Lobatto IIIA 4 &1.26&1.55&1.89&2.32&3.14&4.26\\\hline
\end{tabular}
 \end{table}

\subsection{Borel-Pad\'e-Laplace integrator}
The \ac{BPL} integrator was introduced first in \cite{dina-thesis} and was then used to integrate a certain number of \ac{ODE} and \ac{PDE} \cite{dina-2012}. The integrator was then combined with a finite element solver to solve numerically problems in fluid mechanics \cite{deeb-thesis,deeb:stab-serie} and applied also to preserve geometrical properties of mechanical problems \cite{ahmad_bpl_2014}. Deeb \emph{et al.} \cite{DEEB_2022_bpl} studied this integrator for solving stiff and non-stiff problems. It has been used also for large-time dynamical problems \cite{ahmad_robust_integrators_2019}, and has been implemented in a Proper Generalize Decomposition solver for non-linear diffusion equations \cite{ahmad_pgd_pade}. This integrator has a variant that uses \ac{GFS} \cite{ahmad_icnpaa_2016}. Its efficiency was compared with the \ac{BPL} in \cite{ahmad_comp_bpl_sfg_2015}.

This integrator was inspired by the Borel-Laplace resummation method. Consider first that we are at the instant $\t_n$ and have an approximation $\y_n$. We denote by $Y_{n,0}\coloneqq \y_n$ and for every $\h >0$ we consider that the solution can be written in the form of a time series expansion in the vicinity of $\t_n$ as follows:
\begin{equation}
	\hat \y(\t_n+\h) \coloneqq \sum\limits_{k=0}^{\infty} Y_{n,k} \h^k.
\end{equation}
Approximating the solution using the partial sum may not work every time, because the series could present a small radius of convergence. It could also be divergent, where Gevrey \cite{gevrey_1918} had classified divergences with different orders: we say that a series is divergent of order $r$ if $\lvert Y_{n,k}\lvert  \leqslant CA^k (k!)^{1/r}$. We limit our consideration to the first-order Gevrey series. To elaborate the sum of the divergent series, the inverse of the Laplace transform is applied first on the series $\hat \y$. This is called the Borel transform, which leads to an analytic function inside a disc of convergence in the complex plane and presents singularities around the disc. Then, we prolongate it analytically throughout a semi-line that does not present any singularity, and apply the Laplace transform on this prolongation. If the latter increases not faster than an exponential function at the infinity, its Laplace transform is an analytic function that is Gevrey-asymptotic to the initial series $\hat \y$. For more details, we refer to \cite{deeb-thesis}.
\captionsetup[table]{name=Diagram,skip=1ex, labelfont=bf}
\begin{table*}[ht]\centering
\caption{\ac{BPL} algorithm}
\label{alg:borel_pade_laplace}
  {\sf\small
\begin{tabular}{ccrcc}
$\displaystyle\hat{\y}^\op(\t_n+\h)=\sum_{k=0}^{\op}Y_{n,k}\h^k$&    &$\displaystyle \Phi^{\op}_{\h}(\y_n) =\y_n+\h\sum_{i=1}^{N_G}\mathrm{Pd}^{\op}(\xi_i\cdot\h)\cdot\omega_i$
\\\\
 $\left.\begin{array}{c}\\\text{Borel}\\\\\end{array}\right\downarrow$&
 &$\left\uparrow \begin{array}{c}\\\text{Gauss-Laguerre} \\\\\end{array}\right.$
\\\\
$\displaystyle\big(\B\hat{\y}^\op\big)(\xi)=
\sum_{k=0}^{\op-1}\cfrac{Y_{n,k+1}}{k!}\ \xi^k$
&$\overrightarrow{\hspace{.5cm}\text{Pad\'e}\hspace{.5cm}}$
&$\mathrm{Pd}^\op(\xi) = \cfrac{q_0+q_1\xi+\cdots+q_{\dg q}\xi^{\dg q}}{r_0+r_1\xi+\cdots+r_{\dg r}\xi^{\dg r}},$ \\
&& ${\dg q+r}=\op$
\end{tabular}}
\end{table*}
\captionsetup[table]{name=Table,skip=1ex}

Numerically, the series is truncated up to order $\op$ and the \ac{BPL} algorithm is presented in Diagram \ref{alg:borel_pade_laplace}. After applying the Borel transform of the truncated series, Pad\'e approximants are used to extrapolate the obtained series by the function $\mathrm{Pd}^{\op}(\xi)$. Then, we apply the Laplace transform (the inverse Borel transform) to go back to the physical space. It is approximated using Gau\ss-Laguerre quadrature with $N_G$ Gau\ss\,points $\{\xi_i\,\vert\, i\in\S{1}^{N_G}\}$ and weights $\{\omega_i\,\vert\, i\in\S{1}^{N_G}\}$.
To step forward, the time step $\h$ should be determined such that the error of the approximation does not exceed a predefined user tolerance $\tol$. Yet, there is no a priori error estimate and the valid time step is chosen according to the residual error:
\begin{equation}
 \label{residual_Error}
 \Res\Big( \Phi^{\op}_{\h}(\y_n)\Big) = \frac{\d}{\d\h}\Big(  \Phi^{\op}_{\h}(\y_n)\Big) - \f\Big(\t_n +\h,\Phi^{\op}_{\h}(\y_n) \Big),
\end{equation}
such that the latter does not exceed $\tol$. Practically, we start with a starting value $\h_e$ and evaluate the associated residual error. If $\Res\Big( \Phi^{\op}_{\h_e}(\y_n)\Big)<\tol$, we evaluate again the flow for $\h_e = \h_e\times  C$ ($C\sim 1.1$, increase the step by 10 $\%$ of its value) and its residual. We repeat this loop until reaching an error that is bigger than the defined tolerance.
These evaluations have a high cost of computation, as it requires to evaluate the residual several times. Here, we will use the proposed process to produce a numerical solution using double composition of the \ac{BPL} integrator $\Phi_\h^{\BPL}$, then have an error estimate in the imaginary part $\Im\big(\Psi_\h^{\BPL}\big)$ to use it in the adaptivity technique. The coefficients are defined relative to the truncation order $\op$: $\gamma_1 \coloneqq \frac{1}{2} + \i\,\frac{\sin(\pi/(\op+1))}{1 + \cos(\pi/(\op+1))} $.


\section{Numerical tests: Application to \ac{ODE}s}
\label{sec4}
In this section, we test the imaginary part and its potential to produce numerical solutions with \ac{ATS} using double composition of a basic integrator. To do that, we consider academic \ac{ODE}s where we do have exact solutions or some conserved quantities to compare with numerical ones. The adaptivity of the time step using the imaginary part will be demonstrated, compared to, on one hand, the error estimation if provided by the numerical scheme (as the ERK), and on the other hand to the evolution of the time step if the exact error of the numerical solution were already known.

\subsection{The first example}
\begin{figure}[!ht]
\centering
 \begin{subfigure}{0.75\textwidth}
  \includegraphics[width=\textwidth]{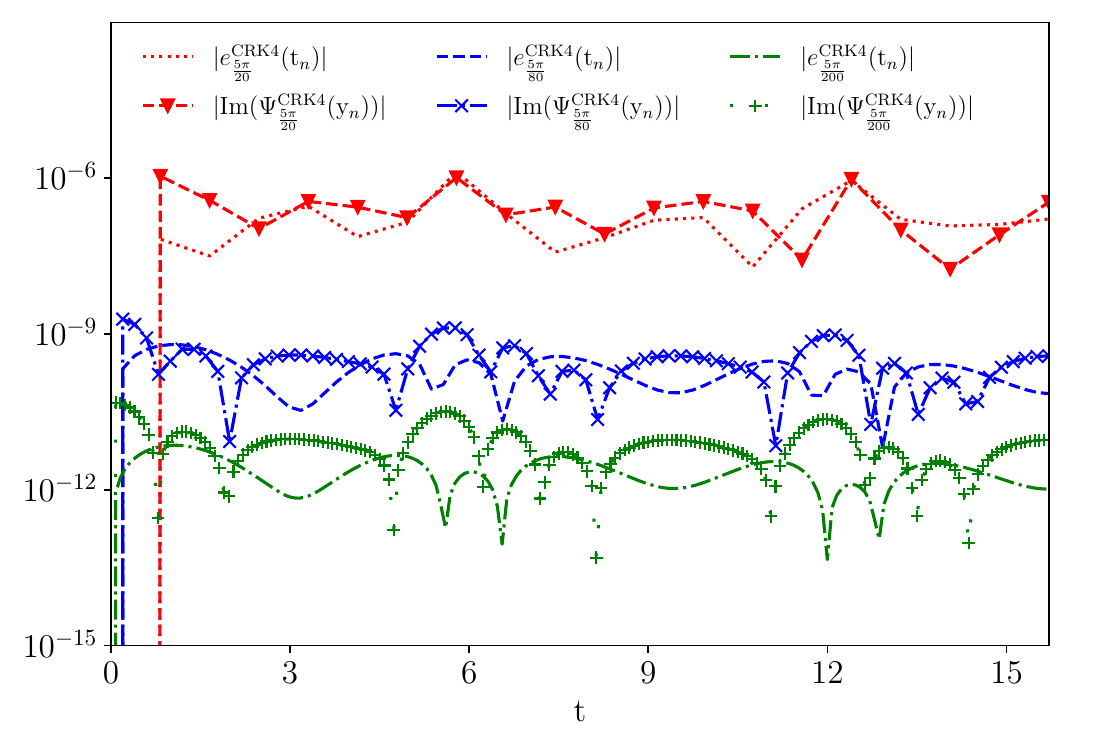}
  \caption{Exact error and imaginary part of the composition of the \ac{RK}4.\label{Fig1_1}}
 \end{subfigure}
 \begin{subfigure}{0.75\textwidth}
    \includegraphics[width=\textwidth]{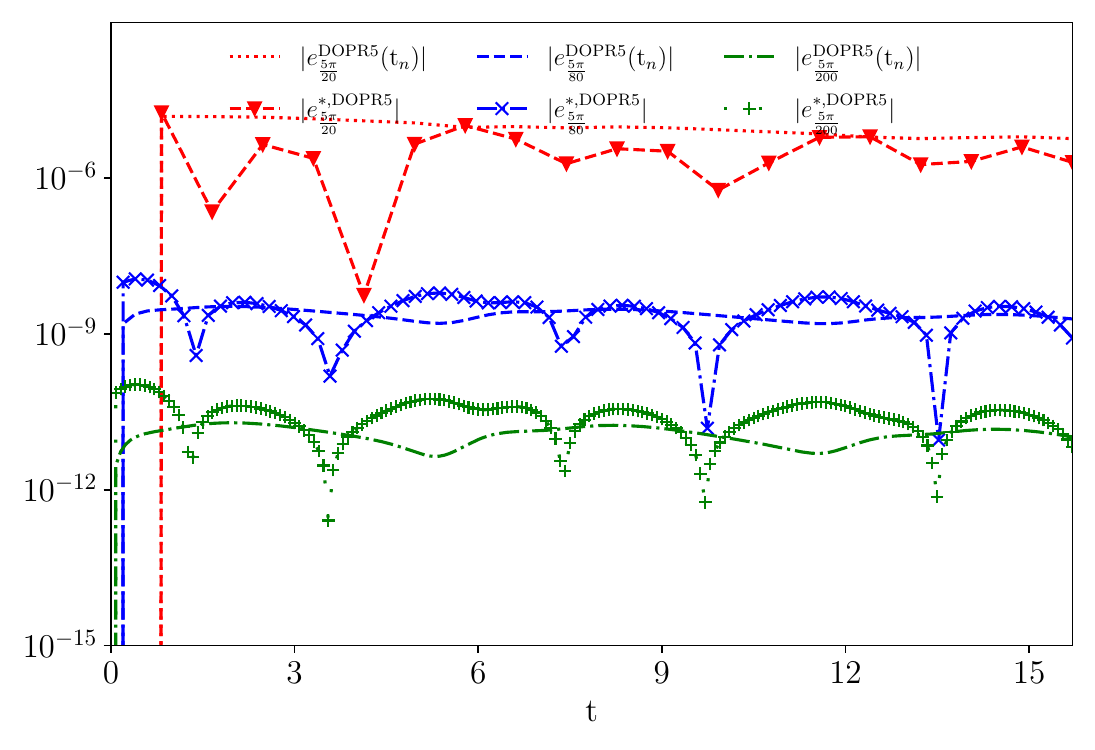}
    \caption{Comparison between the estimation via $\DOP5$ and the exact error. \label{Fig1_2}}
 \end{subfigure}
\caption{Comparison of the exact error obtained by the real part of $\Psi_\h^{\CRK{4}}$ and its imaginary part for different sizes of $\h$ (upper). The lower panel represents the comparison of the error estimation obtained by the Dormand and Prince scheme ($\DOP5$) with the exact error of its approximation. \label{Fig1}}
\end{figure}

We consider first an example of an \ac{IVP} with the initial condition $\y(0) = 1$ and $\f(\t,\y) = e^{-\lambda\y} + \sin(\t)$, where we seek for approximations over the interval $[0,5\pi]$. The exact solution to this equation is given below:
\begin{align*}
	\y(\t) &= \lambda^{-1}\log\Big(\lambda \t  +  \big(\lambda^2 g(\t)  + \e^{\lambda(1+\y_0)}\big) \times \e^{-\lambda \cos(\t)}   \Big), \\
	g(\t) &= \int_0^\t \h \sin(\h)\e^{\lambda \cos(\h)} \d\h.
\end{align*}
We compute a numerical approximation by considering the real part of composing two times the classical fourth-order \ac{RK} scheme, where the imaginary part is compared with the exact error $\e(\t)$ as having the exact solution. \cref{Fig1} shows this comparison for different fixed time steps $\h \in \{{5\pi}/{20}, {5\pi}/{80}, {5\pi}/{200}\}$. We can see that the imaginary part (dashed line with marker) is in the same range of the values of the exact error (dashed line without marker) between the exact solution and the numerical approximation obtained by the real part of the composition:
\begin{equation*}
	 \e^{\CRK{4}}_\h(\t_n) = \y(\t_n) - \Re\big( \Psi^{\CRK{4}}_\h(\y_{n-1}) \big).
\end{equation*}

To check the performance of the imaginary part in error estimation, we compare it with other strategies as the one obtained by Embedded \ac{RK} methods. The scheme $\DOP5$ is employed for its fifth order of accuracy having the same of the composed flow $\Psi_{\h}^{\CRK4}$. The error estimation is of order four. First, we show in \cref{tab01} the CPU time needed to reach the simulation with the associated global error. We see that for all the cases, we reach a higher precision using the composition technique with less CPU time.

\begin{table}[!ht]
\centering
 \caption{Global error versus CPU time for both schemes.\label{tab01}}
 \renewcommand\arraystretch{1.2}
 \begin{tabular}{c|c|ccc}
 $\h$ && $\frac{5\pi}{20}$ &$\frac{5\pi}{80}$ &$\frac{5\pi}{200}$ \\ \hline
  \multirow{2}{*}{$\Phi_\h^{\DOP5}$} &CPU [sec]& $7.41\times10^{-4}$ & $2.9\times10^{-3}$ & $7.15\times10^{-3}$  \\
  & $\overline{e}_\h$ & $4.075\times10^{-5}$ & $9.72\times 10^{-9}$& $5.043\times10^{-11}$ \\ \hline
   \multirow{2}{*}{$\Re\big(\Psi_{\h}^{\CRK4}\big)$}& CPU [sec] &$3.56\times 10^{-4}$& $1.399\times 10^{-3}$& $3.449\times 10^{-3}$ \\
   & $\overline{e}_\h$ & $1.56\times 10^{-6}$ & $1.04\times 10^{-9}$ & $1.18\times 10^{-11}$
 \end{tabular}
\end{table}

For additional details, we present in \cref{tab0} the global ratio between the exact error relative to the approximation associated by every scheme: the real part of $\Psi_\h$ and the estimation by the scheme $\DOP5$, over the estimated ones:
\begin{equation}
 \label{global_ratio}\
 \int_0^{\T{}}\left\|\frac{\e^{\CRK4}(\t)}{\Im\big(\Psi^{\CRK4}_\h)}\right\| \d\t \quad \text{or} \quad
 \int_0^{\T{}}\left\|\frac{\e^{\DOP5}(\t)}{\e^{*,\DOP5}(\t)}\right\| \d\t,
\end{equation}
by the imaginary part $\Im\big(\Psi_\h^{\CRK4}\big)$, on one hand or by $\e^*_n$ provided by the scheme $\DOP5$ on the other hand. The more the ratio is close to one, the more the strategy of estimating the error is better.
\begin{table}[!ht]
\centering
 \caption{Global ratio between the exact error and the estimated one by both schemes.\label{tab0}}
 \renewcommand\arraystretch{1.2}
 \begin{tabular}{c|ccc}
 $\h$ & $\frac{5\pi}{20}$ &$\frac{5\pi}{80}$ &$\frac{5\pi}{200}$ \\\hline
  $\Phi_\h^{\DOP5}$ & 126.462&  6.291 &1.168 \\\hline
  $\Phi_{\h}^{\CRK4}$ & 1.517& 2.418  & 1.419
 \end{tabular}
\end{table}
The performing of the strategy of the imaginary part looks much more better for larger time steps, while it is competing with $\DOP5$ when the time step is smaller.

\subsection{Example with the Lambert function}
Consider a differential equation with $\f(\t,\y) \equiv \y^2 -\y^3$ with the initial condition $\y(0) = \delta$. The solution to this \ac{IVP} is given explicitly by the function $\y(\t) = \cfrac{1}{W(d e^{\,d-\t}) +1}$, where $d \coloneqq {1}/{\delta}-1$ and $z\mapsto W(z)$ is the Lambert function defined as the solution to the implicit equation $We^W = z$. We choose to approximate the solution to this equation over the interval $[0,\frac{2}{\delta}]$, the solution represents a stiff abrupt variation around $\t = {1}/{\delta}$. For this purpose, we select various numerical schemes and compute a numerical solution using the composition technique with a fixed time step to test how the imaginary parts reproduce the error pattern of the numerical solution without applying the adaptivity in \cref{alg-adap}. Two schemes will be composed: \ac{RK}2 and \ac{RK}4. To compare their performance, two schemes are employed: the Bogacki-Shampine that produces an error estimation of order three, denoted by $\Phi_\h^{\BS3}$, to be compared with $\Im\big(\Phi_\h^{\CRK2}\big)$, and the Dormand-Prince scheme, denoted by $\Phi_\h^{\DOP5}$, and producing an error estimation of order four to be compared with $\Im\big(\Phi_\h^{\CRK4}\big)$.

\begin{figure}[!ht]
\centering
\begin{subfigure}{0.49\textwidth}
    \includegraphics[width=\textwidth]{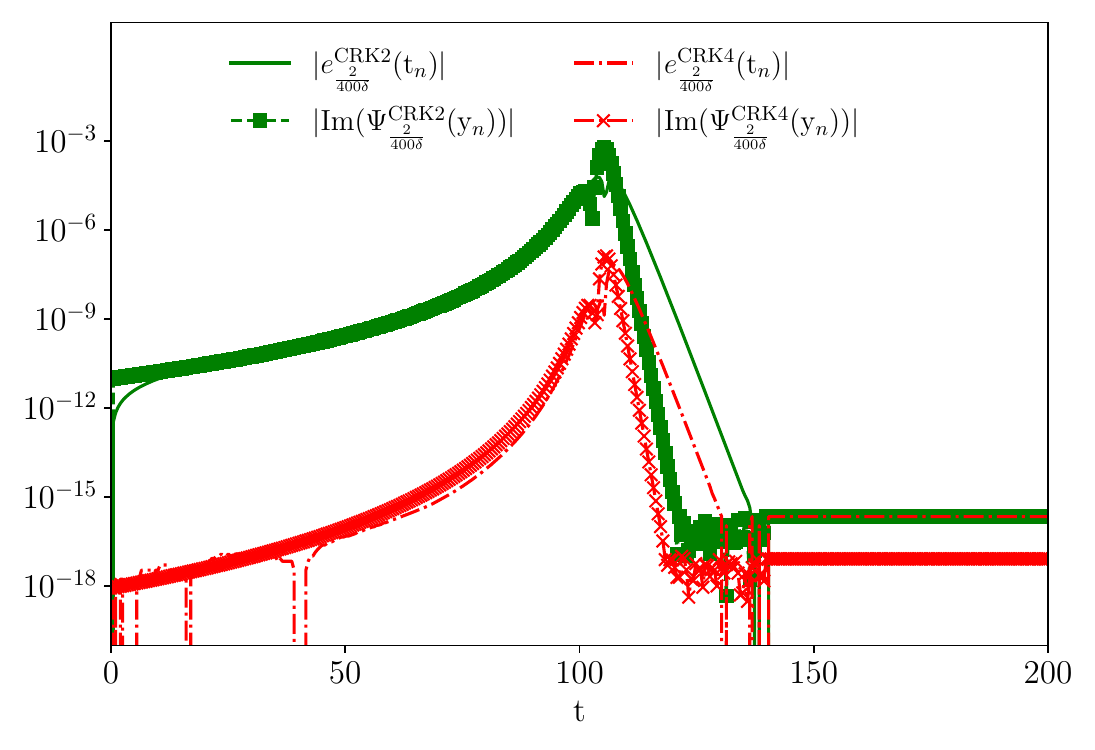}
    \caption{Exact error and imaginary part.}
    \label{Fig21}
\end{subfigure}
\begin{subfigure}{0.49\textwidth}
    \includegraphics[width=\textwidth]{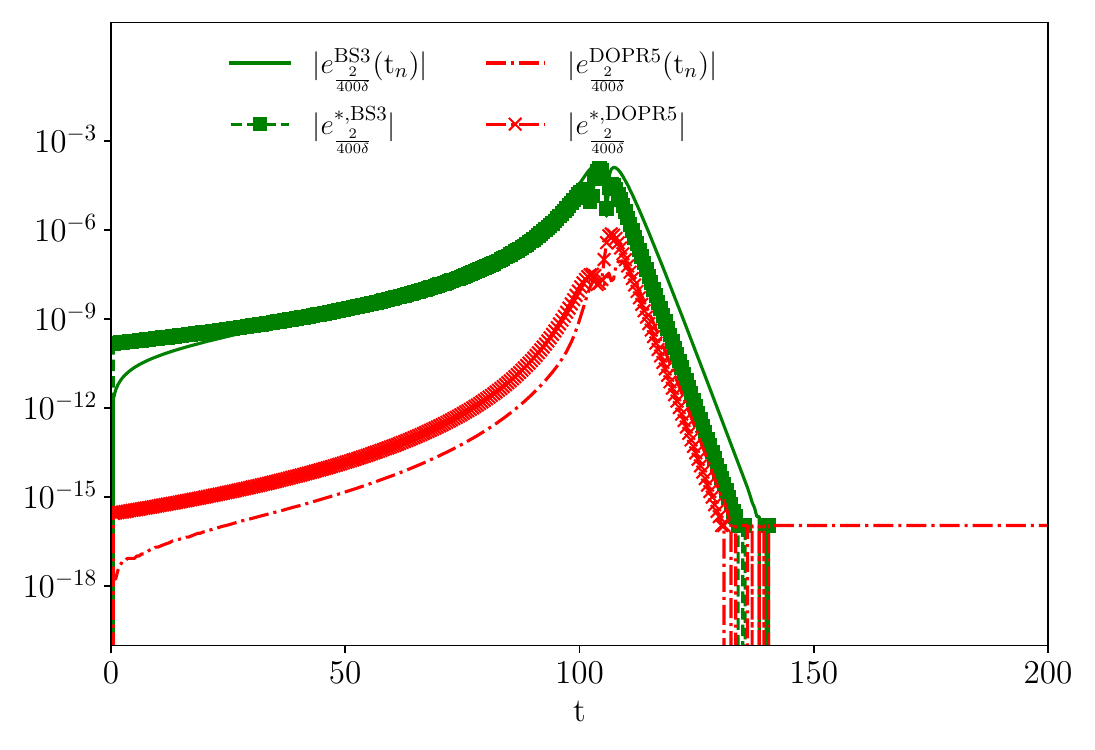}
    \caption{The time step evolution.}
    \label{Fig22}
\end{subfigure}
\caption{Evolution of the numerical solution obtained by a double composition of \ac{RK}2 and \ac{RK}4 for the Lambert problem.}
\label{Fig2}
\end{figure}

First, We present in \cref{Fig2} (left panel) the evolution of the error of the approximation for $\delta = 0.01$, obtained $\Phi_\h^{\CRK2}$ with $\alpha=1/2$ and the famous and $\Phi_\h^{\CRK4}$.
The imaginary part associated with the approximation is also plotted in this figure.
We can see how the imaginary part follows the exact error for both schemes and has the same pattern: it increases in an exponential way to reach a maximum around $\t=1/\delta$, then decreases drastically to reach zero machine precision. This will be used in the time step variation. \cref{Fig22} presents the evolution of the exact error and the estimated ones obtained by both \ac{ERK}: $\BS3$ and $\DOP5$.

To check the performance of the error estimation, the simulation is ran for different time step by the four presented schemes, where the global ratio defined in \cref{global_ratio} is evaluated at every case and results are presented in \cref{tab03}. We can check that the error estimation by the composed schemes present a global ratio closer to one in the most of the cases. In the simulation test, the performance of evert scheme is assessed.  \cref{tab02} presents the global error with the CPU time needed to achieve the simulation using different schemes with different time steps. We conclude also that composed schemes outperform the \ac{ERK} presenting the same orders by achieving the simulation with higher accuracy and less computational time.

\begin{table}[!ht]
\centering
 \caption{Global ratio between the exact error and the estimated one by used schemes.\label{tab03}}
 \renewcommand\arraystretch{1.2}
 \begin{tabular}{c|ccc}
 $\h$ & $\frac{2}{100\delta}$ &$\frac{2}{200\delta}$ &$\frac{2}{400\delta}$ \\\hline
  $\Phi_\h^{\BS3}$  &  0.458 & 0.461 & 0.463 \\
  $\Phi_\h^{\CRK2}$ & 0.467& 0.871& 1.681 \\
  $\Phi_\h^{\DOP5}$ & 0.0223& 0.0278& 0.0393 \\
  $\Phi_\h^{\CRK4}$ &  0.0717& 0.0983& 0.414
  \end{tabular}
  \end{table}

\begin{table}[!ht]
\centering
 \caption{Global error versus CPU time for different schemes.\label{tab02}}
 \renewcommand\arraystretch{1.2}
 \begin{tabular}{c|c|ccc}
 $\h$ && $\frac{2}{100\delta}$ &$\frac{2}{200\delta}$ &$\frac{2}{400\delta}$ \\ \hline
 \multirow{2}{*}{$\Phi_\h^{\BS3}$} &CPU [sec]& $1.98\times10^{-3}$ & $4.07\times10^{-3}$ & $7.57\times10^{-3}$  \\
  & $\overline{e}_\h$ & $2.93\times10^{-2}$ & $2.89\times 10^{-3}$& $3.24\times10^{-4}$ \\ \hline
   \multirow{2}{*}{$\Phi_{\h}^{\CRK2}$}& CPU [sec] &$5.8\times 10^{-4}$& $1.14\times 10^{-3}$& $2.23\times 10^{-3}$ \\
   & $\overline{e}_\h$ & $1.06\times 10^{-2}$ & $1.09\times 10^{-3}$ & $1.24\times 10^{-4}$\\ \hline
 \multirow{2}{*}{$\Phi_\h^{\DOP5}$} &CPU [sec]& $3.04\times10^{-3}$ & $5.81\times10^{-3}$ & $1.11\times10^{-2}$  \\
  & $\overline{e}_\h$ & $8.85\times10^{-4}$ & $1.04\times 10^{-5}$& $1.45\times10^{-7}$ \\ \hline
   \multirow{2}{*}{$\Phi_{\h}^{\CRK4}$}& CPU [sec] &$7.76\times 10^{-4}$& $1.54\times 10^{-3}$& $3.05\times 10^{-3}$ \\
   & $\overline{e}_\h$ & $1.74\times 10^{-4}$ & $3.62\times 10^{-6}$ & $9.72\times 10^{-8}$
 \end{tabular}
\end{table}

\begin{figure}[!ht]
\centering
\begin{subfigure}{0.45\textwidth}
    \includegraphics[width=\textwidth]{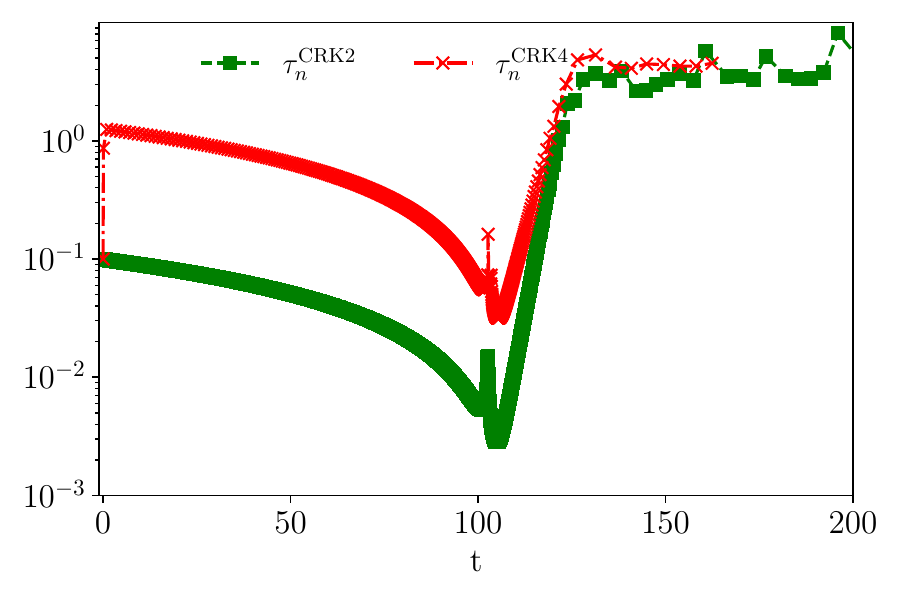}
    \caption{Composition of the second and fourth order \ac{RK}.}
    \label{Fig31}
\end{subfigure}
\begin{subfigure}{0.45\textwidth}
    \includegraphics[width=\textwidth]{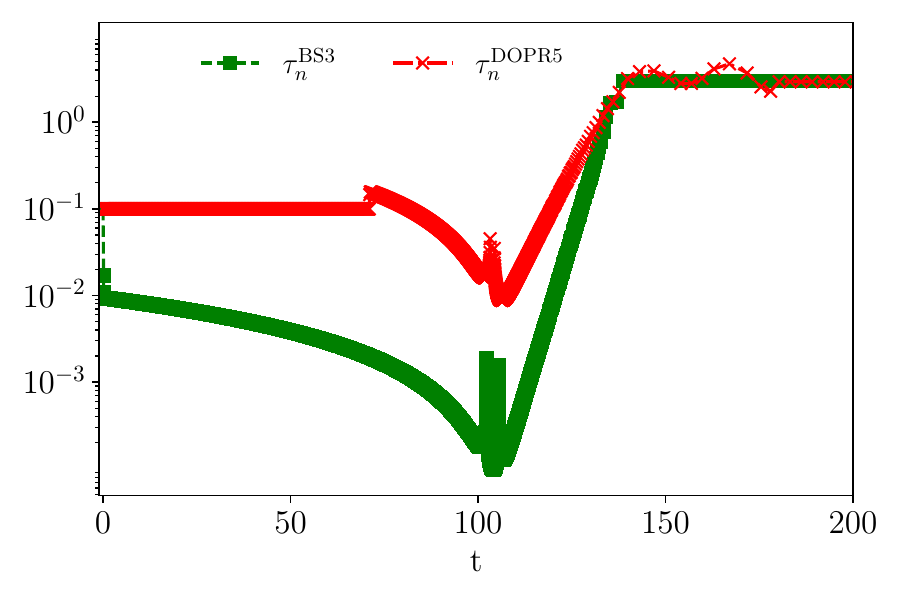}
    \caption{Bogacki-Shampine and Dormand-Prince.\label{Fig32}}
\end{subfigure}
\caption{Evolution of the time steps for the Lambert problem solved by the composition technique (left panel) and by two \ac{ERK} methods: Bogacki-Shampine and Dormand-Prince. $\tol = 10^{-10}$}
\label{Fig3}
\end{figure}

We use now the imaginary part to adapt the time step, where we show in \cref{Fig3} (right panel) how its dynamics within the time for both schemes. The simulation is done also for $\delta = 0.01$ and the user tolerance fixed in the adaptivity formula \eqref{adapt-form} is equal to $\tol = 10^{-10}$ for both schemes. We start with an initial time discretization step $\h_0 = 10^{-1}$ in both cases.

\begin{figure}[!ht]
\centering
\begin{subfigure}{0.45\textwidth}
    \includegraphics[width=\textwidth]{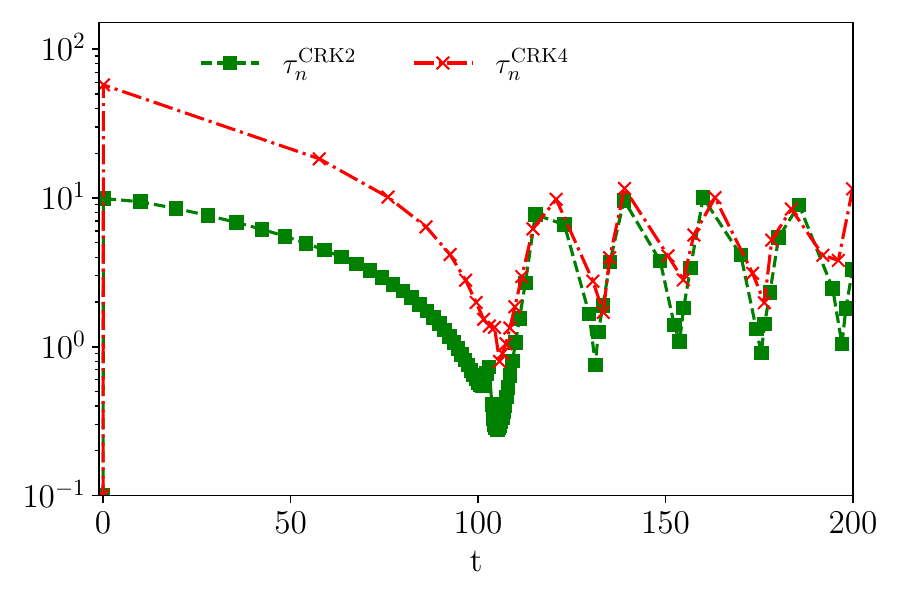}
    \caption{Composition of the second and fourth order \ac{RK}.}
    \label{Fig31_1}
\end{subfigure}
\begin{subfigure}{0.45\textwidth}
    \includegraphics[width=\textwidth]{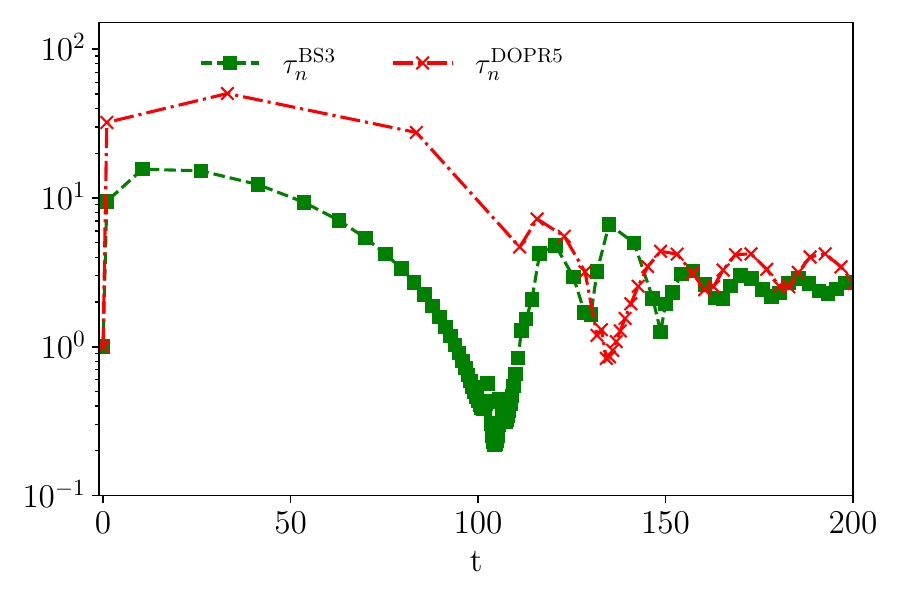}
    \caption{Bogacki-Shampine and Dormand-Prince.\label{Fig32_1}}
\end{subfigure}
\caption{Evolution of the time steps for the Lambert problem. $\tol = 10^{-5}$}
\label{Fig3_1}
\end{figure}

The \ac{ERK} are used here also to adapt the time step by both schemes: $\BS3$ and $\DOP5$ using the same features of the simulation.
For the composition of \ac{RK}4, the time step jumps directly to $\h_1 \approx 1$ and then gets smaller at every iteration before reaching a minimum value at the point $\t\sim1/\delta$. This is encouraging as we observe the same pattern of the evolution of the time step In the case of using \ac{ERK} (see \cref{Fig32}).
Then, both compositions show an increasing time step in a drastic way before oscillating around $\h_n=5$, which demonstrates the utility of the imaginary part as a tool for the \ac{ATS} technique. Results for $\tol = 10^{-5}$ are plotted in \cref{Fig3_1}.

\subsection{Lotka-Volterra problem}
Consider the Lotka-Volterra problem given in the \ac{ODE} system below:
\begin{equation}
\label{lotka-system}
	\left\lbrace
	\begin{array}{l}
	\dot\u = \alpha \u - \beta \u\v,\\
	\dot\v = -\delta \v + \eta \u\v,
	\end{array}
	\right.
	\quad
	\alpha,\beta,\delta,\eta \in \mathds{R}^+.
\end{equation}
This system models the dynamics of two populations: predators ($\u$) and preys ($\v$), where $\alpha$ is the preys production rate, $\delta$ is the predators dying rate, $\beta$ is the rate of decreasing of preys population because of predators, and $\eta$ is the rate of increasing of predators population thanks to preys. The System has the first integral $F(\u,\v)$ given by:
\begin{equation*}
F(\u,\v) = \beta\v + \eta\u  - \alpha\log(\v) - \delta\log(\u).
\end{equation*}
Hence, we have $F(\u(\t),\v(\t)) \equiv F(\u(0),\v(0))$, for all solutions existence time.
\begin{figure}[!ht]
\begin{subfigure}{\textwidth}
    \centering
    \includegraphics[width=0.7\textwidth]{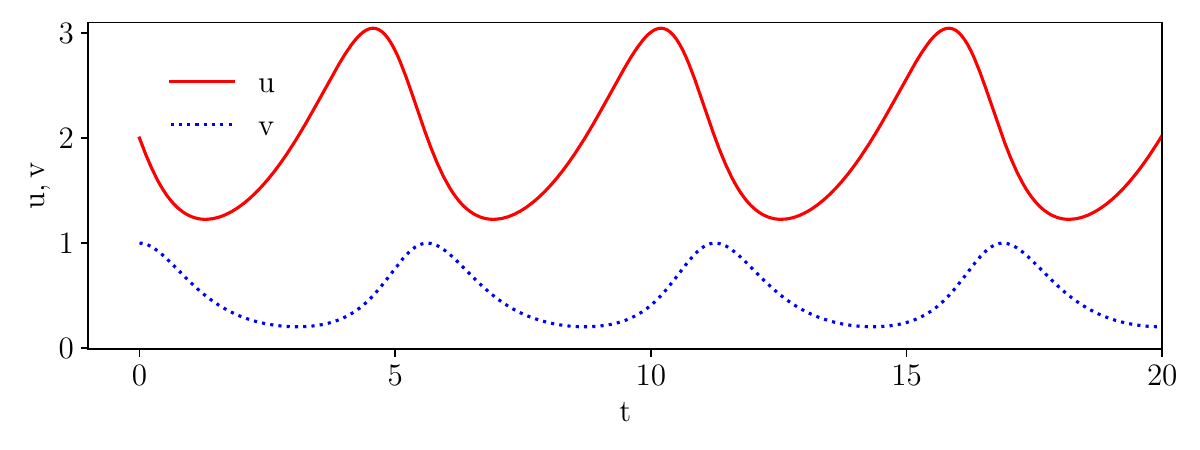}
    \caption{Reference solution.}
    \label{Fig61}
\end{subfigure}
\begin{subfigure}{\textwidth}
    \centering
    \includegraphics[width=0.7\textwidth]{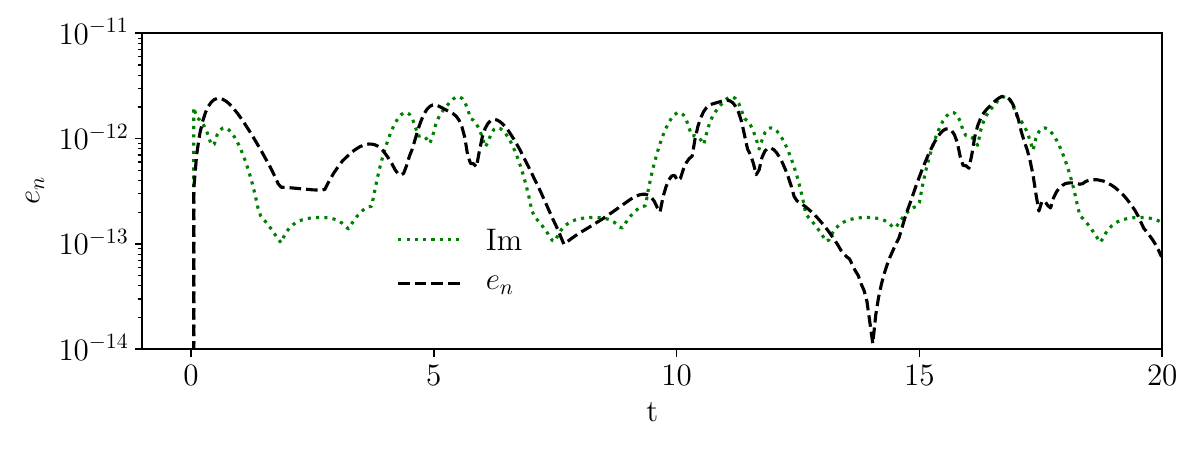}
    \caption{Exact error and imaginary part.}
    \label{Fig62}
\end{subfigure}
\caption{Plot of a reference solution to the Lotka-Volterra System \eqref{lotka-system} (upper panel), the exact error and of the imaginary part of the numerical solution (lower panel) obtained by the composition of the \ac{GRK}2.}
\label{fig:lotka}
\end{figure}
We consider the case where $\u(0) = 2$, $\v(0)=1$ and the solution is sought over the interval $[0,20]$. The simulation, plotted in \cref{fig:lotka}, is done using a fixed time step $\h = 0.5$ and using the composition of the \ac{GRK}2 scheme defined in \cref{sec-GRK2}. We plot in \cref{Fig61} a reference solution obtained with a classical fourth-order \ac{RK} scheme with a time step $\h=10^{-5}$, which is used later to compute the error of the numerical approximation obtained with the real part of the composition of \ac{GRK}2. This error is plotted in \cref{Fig62} and compared with the imaginary part error estimation of the composed numerical flow. We can see how the latter lies in the same range of the error and its pattern resembles closely to that of the exact error.

%

\subsection{Duffing-Van der Pol equation}
Consider the following second order \ac{ODE}:
\begin{equation*}
 \ddot \u + (\r + \g\u^2)\dot\u + \a \u + \b \u^3 = \c \cos(w\t),
\end{equation*}
with $\r,\a,\b,\c,\g$ and $\w$ are real positive constants. When $\r,\a>0$ and $\g,\b=0$, the equation models a single oscillator. When $0 \notin \{\r,\a,\b\} $ and $\g=0$, it produces the Duffing equation that models an oscillator with damping force, and when $ \r<0$, $\a,\g>0$ and $\b=0$ it will be the Van der Pol oscillator. To apply the proposed strategy of constructing numerical solutions with adaptive time steps, we write the equation as the system of the first-order \ac{ODE}s by considering a new variable $\v = \dot\u$. Thus, with $\y \coloneqq (\u,\v)^\top$, the vector function $\f(\t,\y) \coloneqq \big(\v,-(\r+\g\u^2) \v - \a \u - \b \u^3 + \c\cos(\w \t)\big)^\top$ will represent the system. The integrability of this equation has been studied for some cases of parameters to obtain the qualitative behaviour \cite{HOL-1980,Udw-2014} of solutions. However, we are interested in this study in computing the numerical solution within finite time.

\subsubsection{Van der Pol oscillator}
When $\r, \g=0$, the equation can be written in the Hamiltonian form with $\H(\u,\v) \coloneqq \frac{1}{2}\v^2 + \frac{\a}{2}\u^2 + \frac{\b}{4}\u^4$. The quantity $\H\big(\u(\t),\v(\t)\big)$ is conserved during time evolution for every \ac{IVP} and is equal to $\H\big(\u(0),\v(0)\big)$.
\begin{figure}[!ht]
\begin{subfigure}{\textwidth}
    \centering
    \includegraphics[width=0.7\textwidth]{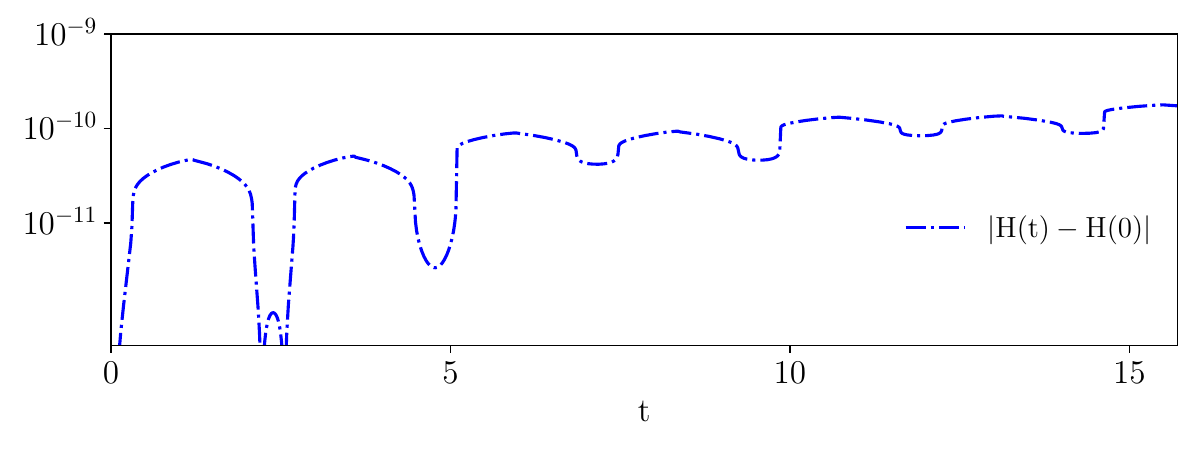}
    \caption{The Hamiltonian error.}
    \label{Fig51}
\end{subfigure}
\begin{subfigure}{\textwidth}
    \centering
    \includegraphics[width=0.7\textwidth]{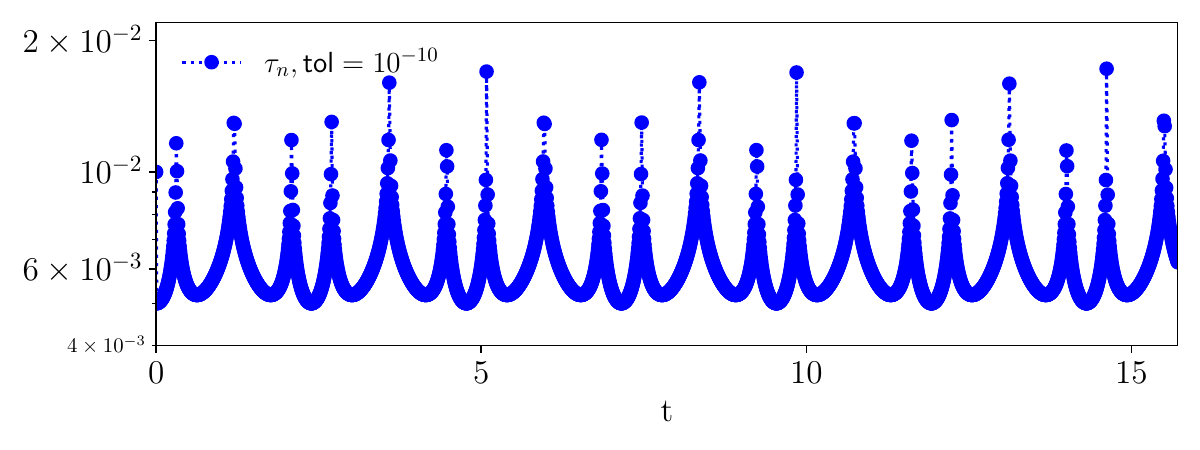}
    \caption{The time step adaptivity.}
    \label{Fig52}
\end{subfigure}
\caption{Evolution of the numerical solution, to the Van der Pol oscillator, obtained by a double composition of \ac{BPL}5 with $\tol = 10^{-10}$.}
\label{Fig5}
\end{figure}
\begin{figure}[!ht]
\begin{subfigure}{\textwidth}
    \centering
    \includegraphics[width=0.7\textwidth]{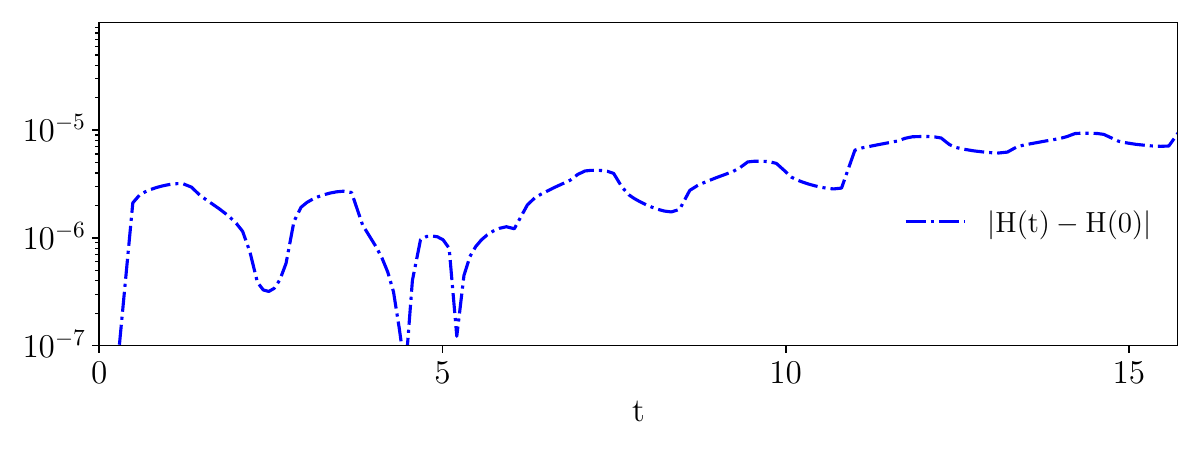}
    \caption{The Hamiltonian error.}
    \label{Fig53}
\end{subfigure}
\begin{subfigure}{\textwidth}
    \centering
    \includegraphics[width=0.7\textwidth]{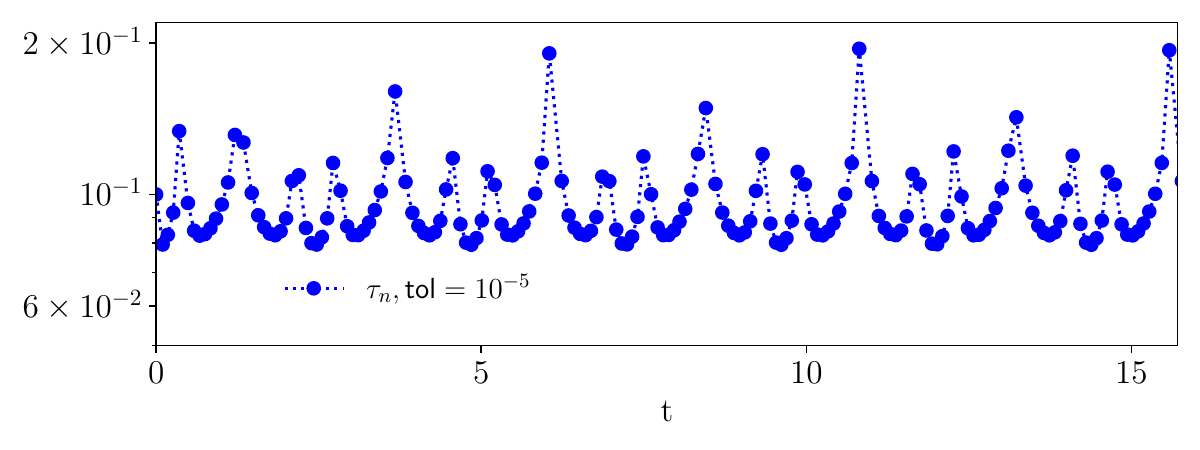}
    \caption{The time step adaptivity.}
    \label{Fig54}
\end{subfigure}
\caption{Evolution of the numerical solution, to the Van der Pol oscillator, obtained by a double composition of \ac{BPL}5 with $\tol = 10^{-5}$.}
\label{Fig55}
\end{figure}

Thus, for numerical simulations, this quantity is used actually to assess the accuracy of the numerical solution and adapt the time step to keep the Hamiltonian conserved up to a given tolerance $\tol$. The imaginary part of the composed BPL integrator of order $N=5$ (BPL5) is used as an error estimate to adapt the time step. The tolerance is set to $\tol = 10^{-10}$, where results are printed in Figure \ref{Fig5}. \cref{Fig51} presents the error of the Hamiltonian obtained via the \ac{BPL}5 and \cref{Fig52} presents the time step evolution regarding its adaptivity based on the imaginary part of the composed flow. We can see that the error of the Hamiltonian is stable within evolution in time, however the token values of the time step span in the range of $[0.004,0.02]$. Another simulation is done with a lower tolerance: $\tol = 10^{-5}$, where results are plotted in \cref{Fig55}.

\subsubsection{Duffing problem}
\begin{figure}[!ht]
\centerline{
\includegraphics[height=3.3cm]{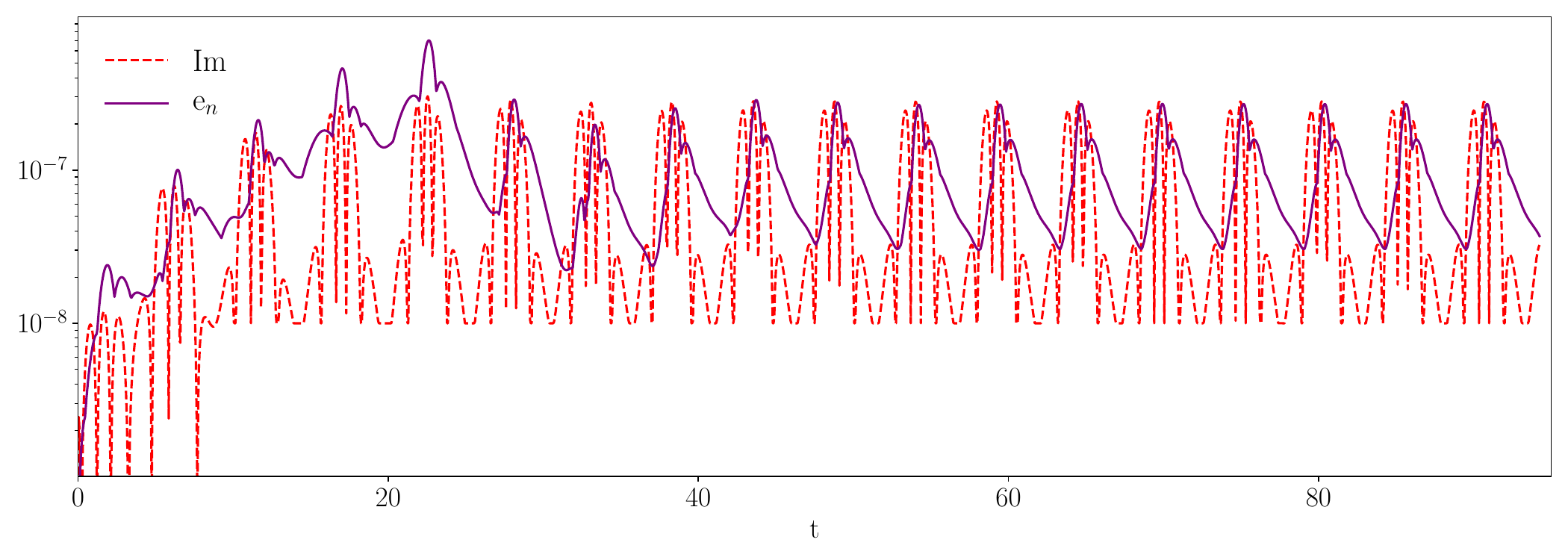}
\includegraphics[height=3.3cm]{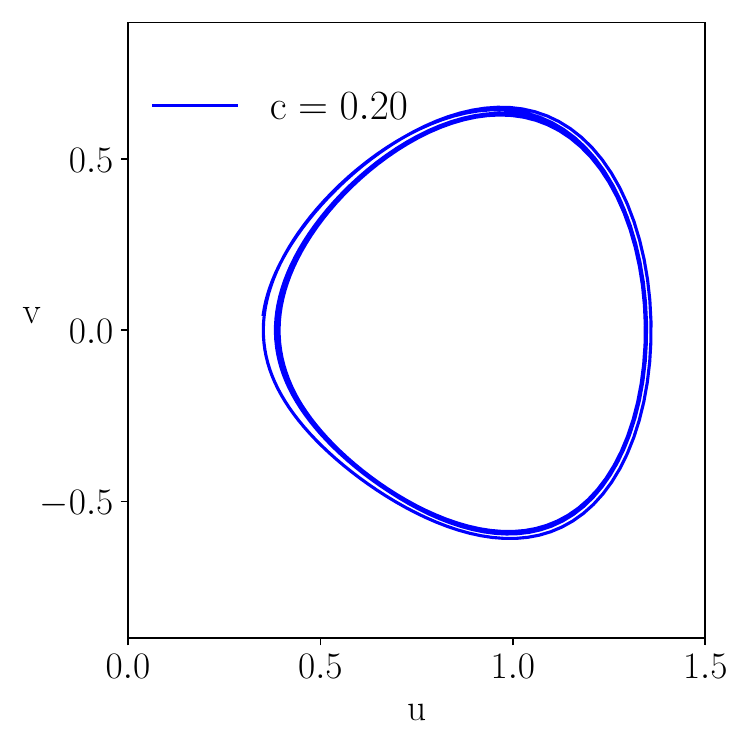}
}
\centerline{
\includegraphics[height=3.3cm]{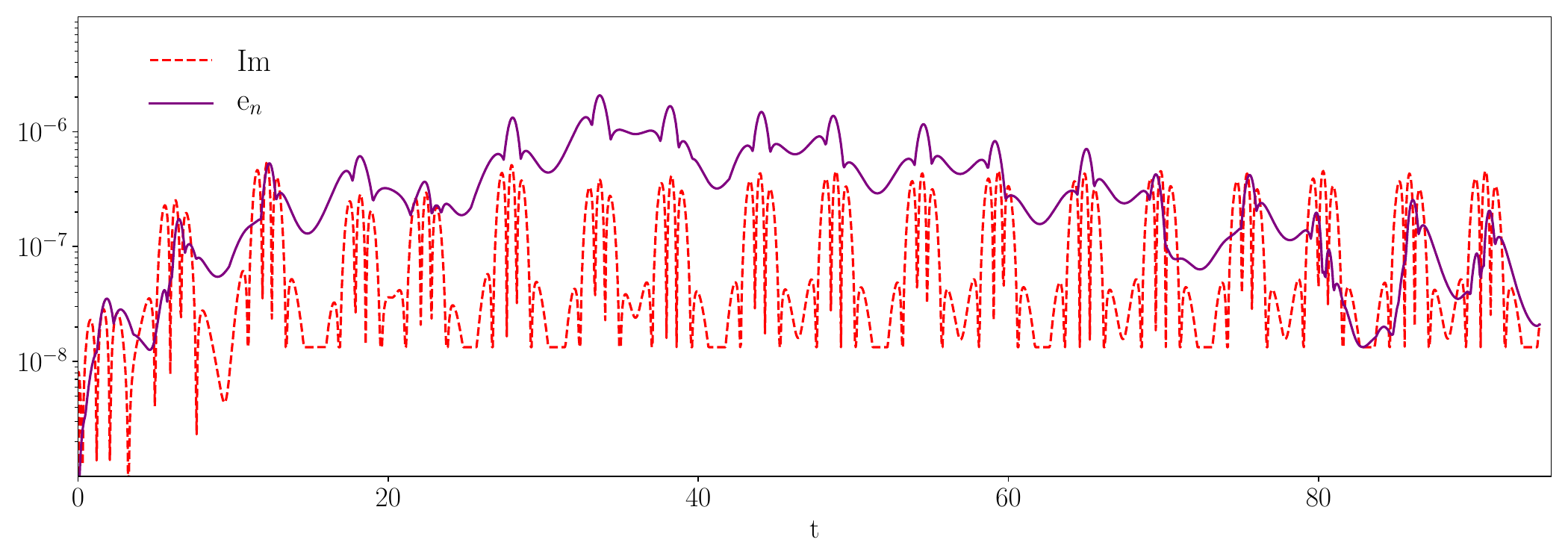}
\includegraphics[height=3.3cm]{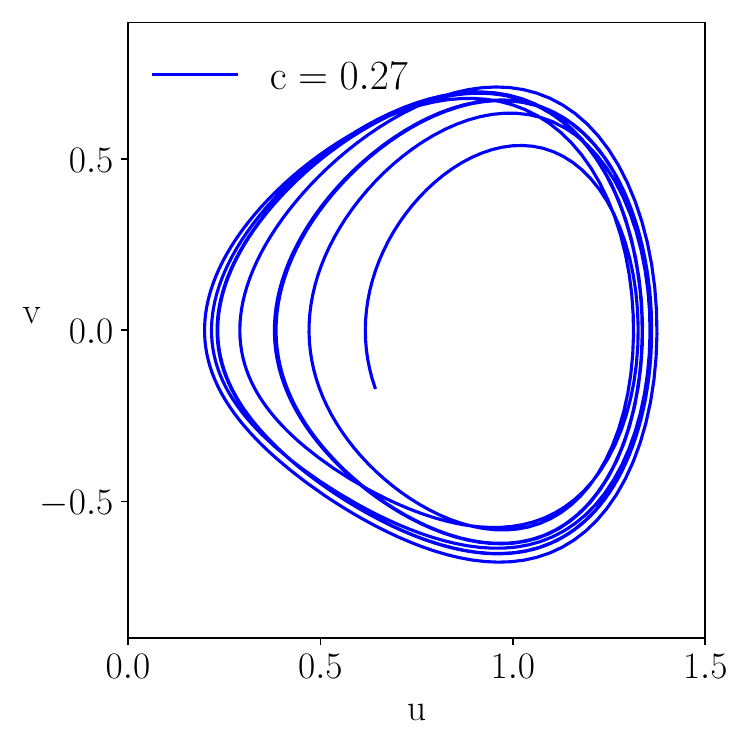}
}
\centerline{
\includegraphics[height=3.3cm]{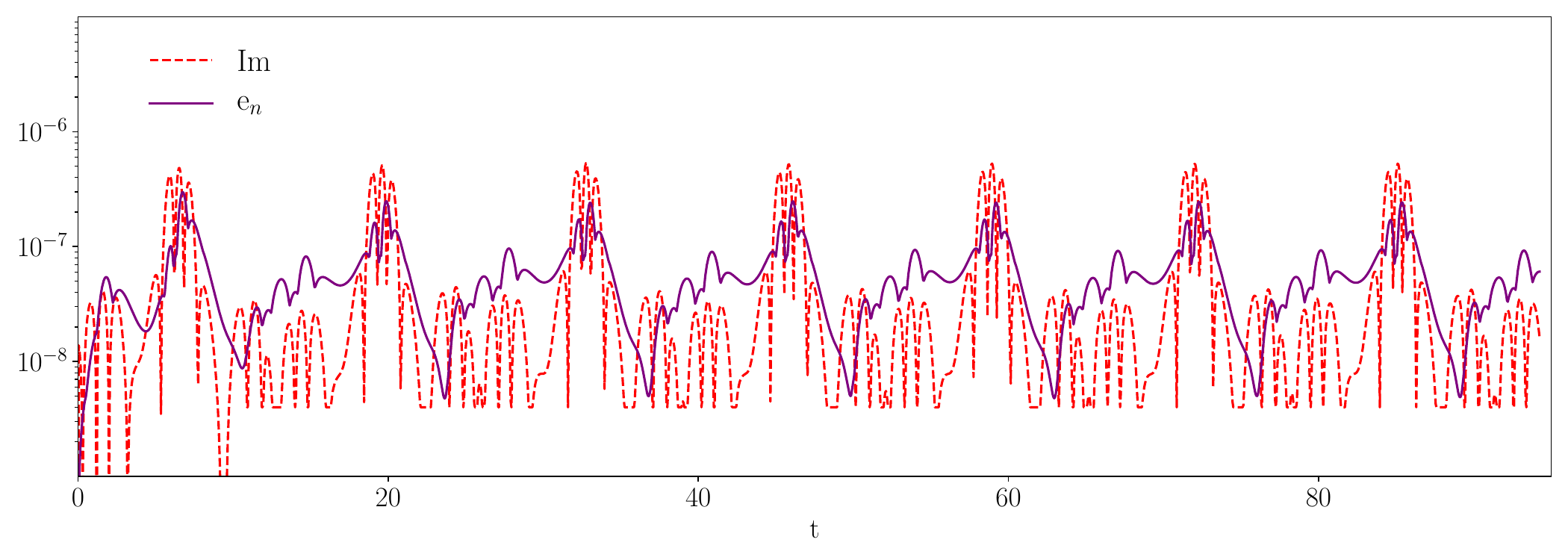}
\includegraphics[height=3.3cm]{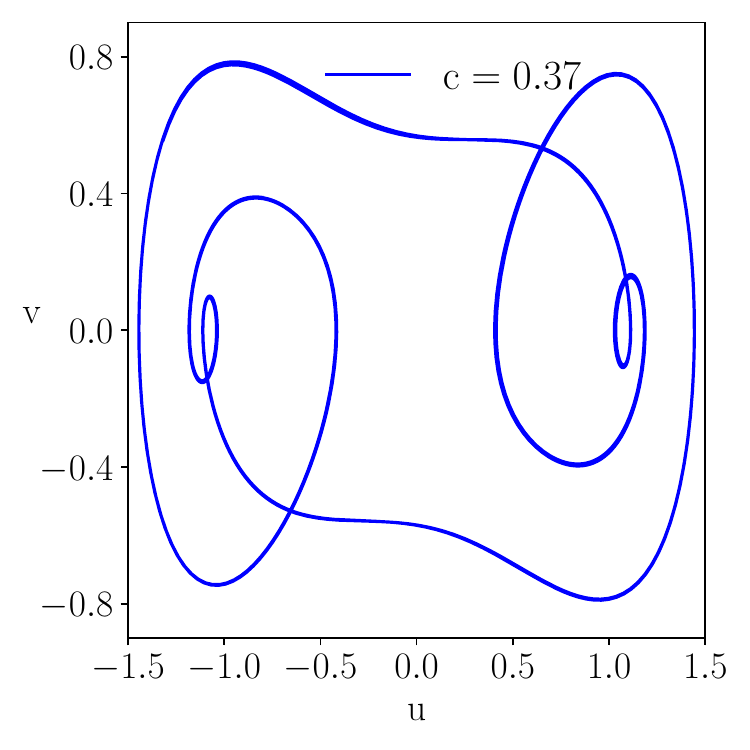}
}
\centerline{
\includegraphics[height=3.3cm]{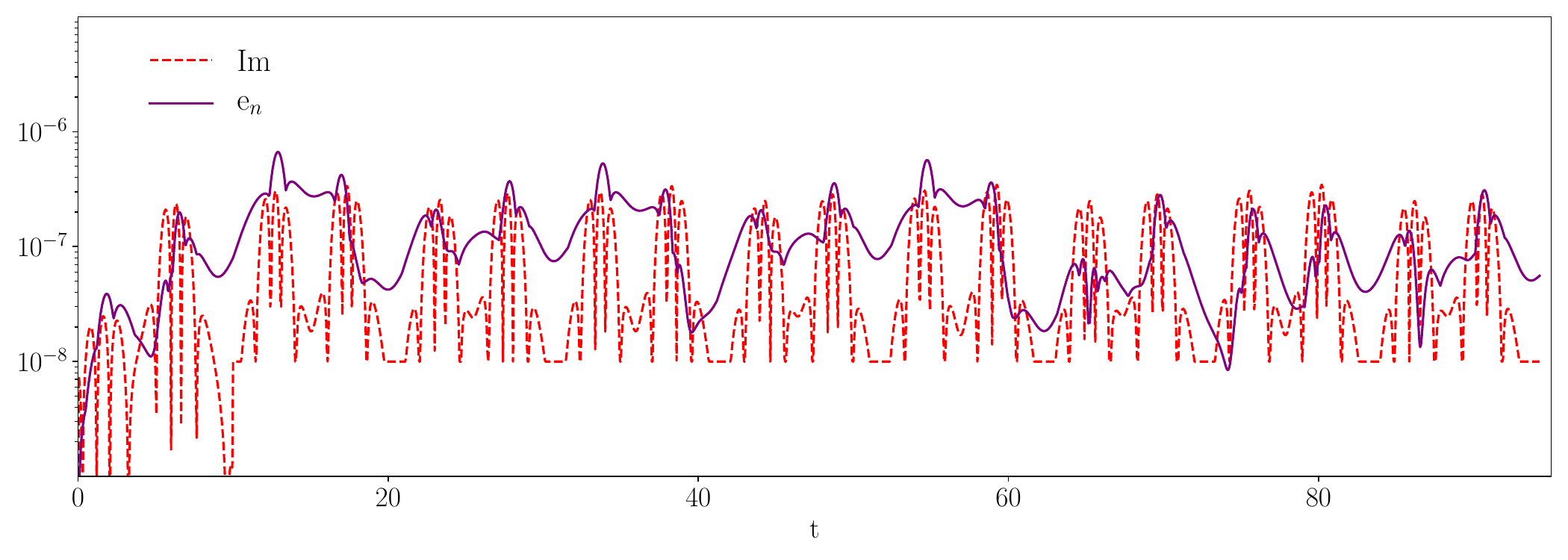}
\includegraphics[height=3.3cm]{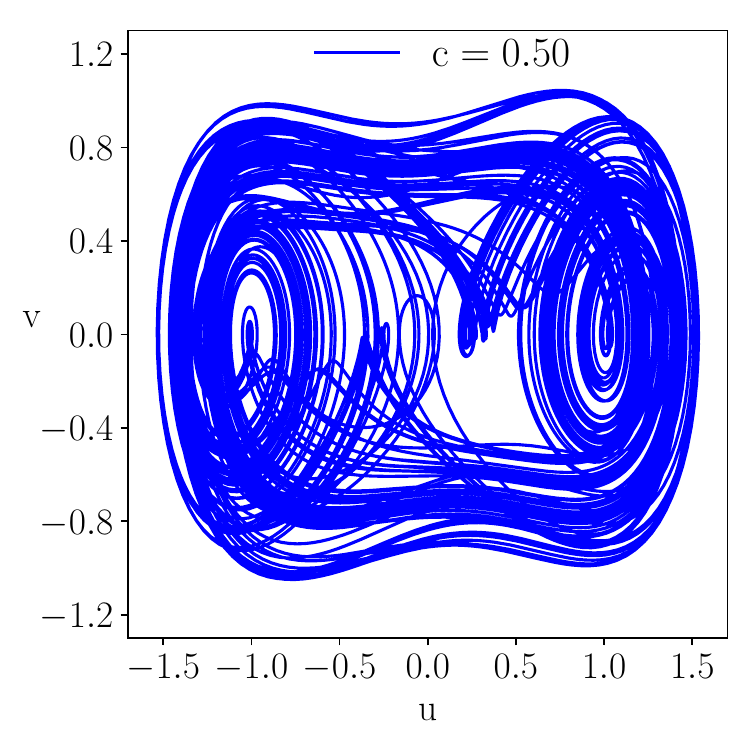}
}
\centerline{
\includegraphics[height=3.3cm]{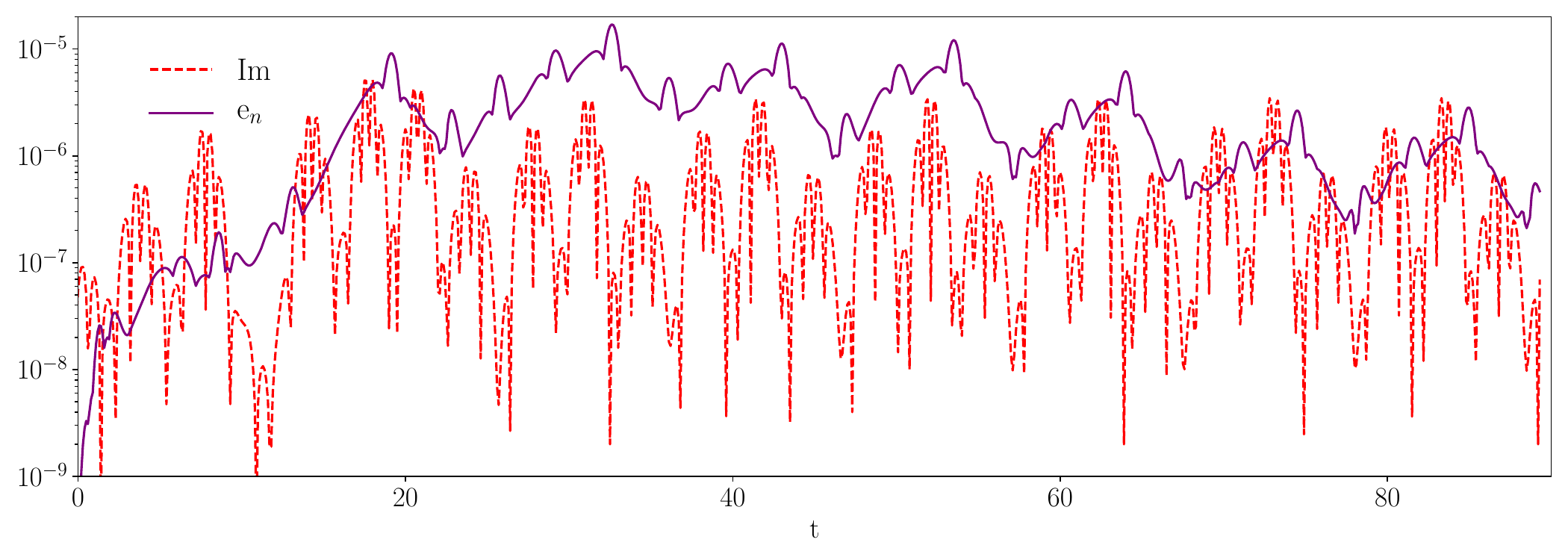}
\includegraphics[height=3.3cm]{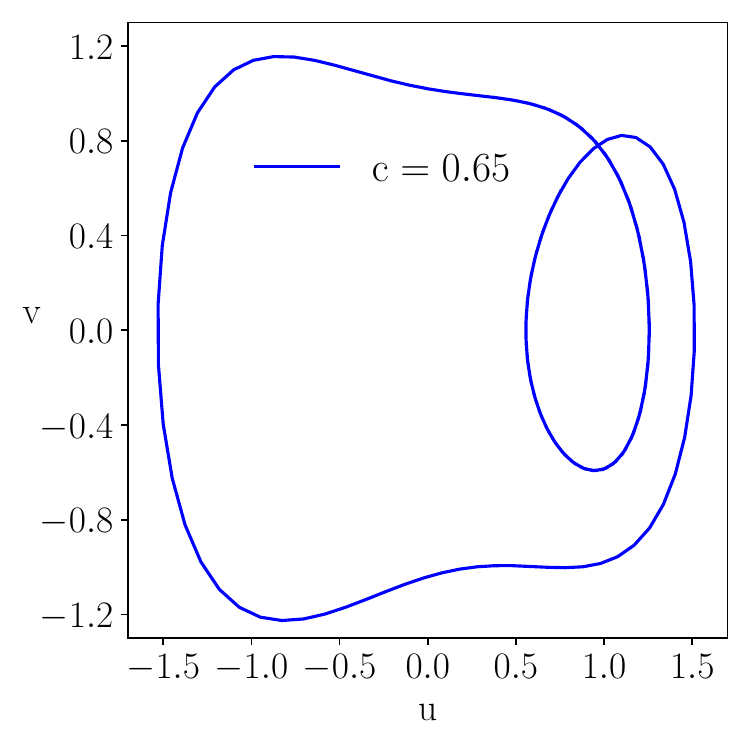}
}
\caption{Duffing problem. Left panel: Evolution of the imaginary part and the error of the approximation obtained by the double composition of \ac{BPL} algorithm with $\op=5$. Right panel: Evolution of the solution in the phase space.}
\label{Fig8-9}
\end{figure}

In this experiment, we study the composition of the \ac{BPL} scheme to integrate the solution of the Duffing problem, where the following parameters are fixed as $\r=0.3,\g=0,\a=-1,\b=1, \w = 1.2$, while $\c$ will be picked from the range $[0.2,0.5]$. The time interval of the simulation is $]0,100[$. The solution of the problem for $\c=0.5$ is chaotic as shown below. The \ac{BPL} solver is processed for $N=5$. Figure \ref{Fig8-9} shows on the right panel the plot, in phase space, of the numerical reference solution obtained by a \ac{GRK}4 scheme with the time step $h=10^{-3}$ for different values of $\c \in\{0.2,0.27,0.37,0.5,0.65\}$.
We use the composition of the \ac{BPL}5 integrator to approximate the solution by its real part, where the error with the reference solution obtained by \ac{GRK}4 is plotted on the left panel with its imaginary part. Graphical results show that the imaginary parts, for cases $\c\in\{0.2, 0.27, 0.37\}$ oscillate, but stay in the same range of the exact error and follow it in a global pattern.

\section{Conclusions and perspectives}
\label{sec5}

\paragraph{Conclusions}
This paper has introduced a new methodology for error estimation in numerical simulations of dynamical systems. By leveraging the complex composition of one-step numerical methods of order $\op$, we have established that the real parts of the resulting compositions furnish approximations of heightened order
$\op+1$, while the imaginary parts deliver reliable error estimates. This approach diverges from traditional practices by utilizing the complex plane, a novelty substantiated by rigorous proof within this work.

The linear stability of the composed schemes showed larger domains of stability in the complex plane, promising an increase in stability of the resulted schemes. On the other hand, the composition technique allows a faster computation, outperforming the basic integrators in providing higher accuracy of numerical approximation with lower time computation.

The practical efficacy of this technique has been thoroughly demonstrated through diverse numerical experiments. The application to several \ac{ODE}s has evidenced the imaginary part's aptitude in mirroring the exact error's behaviour, thereby affirming its theoretical underpinnings. Moreover, the use of the imaginary part in adaptive time-stepping has shown promise, notably in scenarios where conventional error estimates are unattainable or unreliable.

These findings suggest a significant stride forward in numerical analysis, particularly in the adaptive integration of \ac{ODE}s. Future research may explore the extension of these principles to \ac{PDE}s and the potential integration of this method into existing numerical software packages, thereby broadening its applicability to a wider array of scientific computations.

\paragraph{Perspectives}
The methodology presented herein opens multiple promising avenues for future research and application. The intrinsic capacity of the complex composition approach to furnish high-order error estimates in the absence of conventional estimators holds particular promise for the advancement of numerical analysis in computationally intensive fields such as fluid dynamics, climate modeling, astrophysics and quantum mechanics. In the latter a new approach needs to be developed as complex-valued solutions appear (Schrödinger equation), where the theory of multidimensional complex variables is employed.
Further exploration into the application of this method to stiff and multi-scale problems could significantly enhance the robustness and efficiency of simulations in these areas. Moreover, the integration of this error estimation technique with machine learning algorithms may yield adaptive schemes capable of autonomously refining their accuracy in real-time, a frontier that melds traditional numerical methods with modern computational intelligence. Continued development and dissemination of this technique, possibly through open-source numerical libraries, could democratize access to high-precision computational tools, fostering innovation and interdisciplinary collaboration. Ultimately, the theoretical insights gained through this work may also catalyze advancements in the mathematical foundations of numerical error estimation, potentially leading to new theoretical paradigms that align closely with practical computational demands.

\section*{Acknowledgments}

This publication is based upon work supported by the Khalifa University of Science and Technology under Award No. FSU-2023-014.

\section*{List of abbreviations}
\begin{acronym}[]
\acro{ODE}{Ordinary Differential Equation}
\acro{PDE}{Partial Differential Equation}
\acro{IVP}{Initial Value Problem}
\acro{CP}{Cauchy Problem}
\acro{LMS}{Linear Multi-Step}
\acro{BDF}{Backward Difference formula}
\acro{BDFk}{\ac{BDF} of order $k$}
\acro{RK}{Runge-Kutta}
\acro{CRK}{Composed Runge-Kutta}
\acro{ERK}{Embedded-Runge-Kutta}
\acro{IRK}{Implicit \ac{RK}}
\acro{GRK}{Gau\ss-Runge-Kutta}
\acro{ETD}{Exponential-Time Difference}
\acro{DSR}{Divergent Series Resummation}
\acro{GFS}{Generalized Factorial Series}
\acro{BPL}{Borel-Pad\'e-Laplace}
\acro{ATS}{Adaptive Time Stepping}
\end{acronym}

\nomenclature{$\t$}{Time variable}
\nomenclature{$\T{}$}{the final time of simulation}
\nomenclature{$\y(\t)$}{The unknown variable}
\nomenclature{$\f$}{Right hand side of the differential system}
\nomenclature{$\f^{\prime}$}{First total derivative of $\f$}
\nomenclature{$\e(\t)$}{Local error}

\nomenclature{$\phi_t$}{Exact flow of the differential system}
\nomenclature{$\t_n$}{$n$\up{th} instant time}
\nomenclature{$\h_n$}{$n$\up{th} time step}
\nomenclature{$\y_n$}{$n$\up{th} approximation of $\y(\t_n)$}
\nomenclature{$\e_n$}{$n$\up{th} error estimate of $\e(\t_n)$}
\nomenclature{$\overline{\e}_\h$}{Global error of $\e_n$ over the time interval decomposed into sub-interval $\h$}
\nomenclature{$\S{i}^j$}{The set of positive integers numbers between $i$ and $j$ including}
\nomenclature{$\tol$}{User tolerance}
\nomenclature{$\Phi_{\h_n}$}{Numerical flow with a time step $\h_n$}

\nomenclature{$\gamma_1$}{Complex number of the first step integration in the composition}
\nomenclature{$\gamma_2$}{Complex number of the second step integration in the composition}
\nomenclature{$\Re(\y_n)$}{The real part of $\y_n$}
\nomenclature{$\Im(\y_n)$}{The imaginary part of $\y_n$}
\nomenclature{$\Psi_{\h_n}$}{The numerical flow of the double composition of $\Phi_{}$}
\printnomenclature


\appendix
\section{Proof of \cref{Im_12}}
\label{app1}
We have shown first that the coefficients verifying \cref{eq-1} in addition to $\gamma_1 + \gamma_2 = 1$ are given by the following relation:
\begin{equation*}
 \gamma_1 = \frac{1}{2}\Big( 1 + \i \cdot b\Big), \quad \gamma_2 = \overline{\gamma_1}, \quad b= \frac{\sin(\frac{\pi}{\op+1})}{1+\cos(\frac{\pi}{\op+1})}.
\end{equation*}
We first show, by using the binomial formula, that:
\begin{eqnarray*}
\Re{\big(\gamma_1^{p+1}\big)} &=&  \frac{1}{2^{\op+1}}\sum\limits_{n=0}^{[\frac{\op+1}{2}]}\binom{\op+1}{2n}\cdot (\i)^{2n}\cdot b^{2n  },\\
&=&  \frac{1}{2^{\op+1}}\sum\limits_{n=0}^{[\frac{\op+1}{2}]}\binom{\op+1}{2n}\cdot (-\i)^{2n}\cdot b^{2n  },\\
&=& \Re\big(\gamma_2^{\op+1}\big),
\end{eqnarray*}
thus, $\Re\big( \gamma_1^{\op+1}\big) = 0$ using the fact that $\Re\big(\gamma_1^{\op+1} + \gamma_2^{\op+1} \big) =0$.
We multiply now \cref{eq-1} first by $\gamma_1$ and add the result to the multiplication of \cref{eq-1} by $\gamma_2$:
\begin{equation}
\label{eq12}
 \gamma_1^{\op+2} + \gamma_1\gamma_2^{\op+1} + \gamma_2\gamma_1^{\op+1} + \gamma_2^{\op+2}=0.
\end{equation}
Using the theory of complex numbers, we show that:
\begin{eqnarray*}
 \Im\big(\gamma_1\gamma_2^{\op+1}\big) &=& \Im\big(\gamma_1\big)\Re\big(\gamma_2^{\op+1}\big) +\Re\big(\gamma_1\big)\Im\big(\gamma_2^{\op+1}\big)\\
 \Im\big(\gamma_2\gamma_1^{\op+1}\big) &=& \Im\big(\gamma_2\big)\Re\big(\gamma_1^{\op+1}\big) +\Re\big(\gamma_2\big)\Im\big(\gamma_1^{\op+1}\big)
\end{eqnarray*}
By reusing $\Re\big(\gamma_j^{\op+1}\big) = 0$ for $j=1,2$ and that $\Re(\gamma_2) = \Re(\overline{\gamma_1}) = \Re\big(\gamma_1\big)$, we add both equations and simplify the sum to show that:
\begin{equation*}
 \Im\big(\gamma_1\gamma_2^{\op+1} +\gamma_2\gamma_1^{\op+1} \big) = \Re\big( \gamma_1\big)\cdot \Im\big(\gamma_1^{\op+1} + \gamma_2^{\op+1}\big) = 0,
\end{equation*}
which conclude the proof as we have from \cref{eq12}:
\begin{equation*}
 \Im\big(\gamma_1^{\op+2} +  \gamma_2^{\op+2}\big) = -\Im\big(\gamma_1\gamma_2^{\op+1} + \gamma_2\gamma_1^{\op+1} \big)=0.
\end{equation*}


\bibliographystyle{elsarticle-num}
\bibliography{references1}

\end{document}